\documentclass[final,onefignum,onetabnum]{siamart220329}

\usepackage{placeins}
\usepackage{graphicx}

\usepackage{amsmath,amsfonts,amssymb}

\usepackage[pagewise]{lineno}
\usepackage{color}
\definecolor{blue}{rgb}{0.,0.,1.}
\definecolor{red}{rgb}{1.0,0.,0.}
\definecolor{green}{rgb}{0.1,0.5,0.1}

\usepackage[capitalize,nameinlink]{cleveref}
\Crefname{figure}{Figure}{Figures}

\crefformat{equation}{#2(#1)#3}
\crefrangeformat{equation}{(#3#1#4--#5#2#6)}
\crefmultiformat{equation}{#2(#1)#3}{ and #2(#1)#3}
{, #2(#1)#3}{, and #2(#1)#3}
\crefrangemultiformat{equation}{(#3#1#4--#5#2#6)}%
{ and (#3#1#4--#5#2#6)}{, (#3#1#4--#5#2#6)}{, and (#3#1#4--#5#2#6)}

\Crefformat{equation}{#2Equation~(#1)#3}
\Crefrangeformat{equation}{Equations~(#3#1#4--#5#2#6)}
\Crefmultiformat{equation}{Equations~(#2#1#3)}{ and (#2#1#3)}
{, (#2#1#3)}{, and (#2#1#3)}
\Crefrangemultiformat{equation}{Equations~(#3#1#4--#5#2#6)}%
{ and (#3#1#4--#5#2#6)}{, (#3#1#4--#5#2#6)}{, and (#3#1#4--#5#2#6)}

\newcommand{\add}[1]{{\color{blue}#1}}
\newcommand{\remove}[1]{{\color{red} [[#1]] }}

\renewcommand{\add}[1]{#1}
\renewcommand{\remove}[1]{}

\usepackage{color}

\newcommand{\ignore}[1]{}
\newcommand{\eqn}[1]{(\ref{#1})}
\newcommand{\If}{\mathcal{I}_f}
\newcommand{\half}{\frac{1}{2}}
\newcommand{\dt}{\Delta t}
\newenvironment{mat}{\left[ \begin{array}{cccc}}{\end{array}\right]}
\newcommand\bcm{\begin{mat}}
\newcommand\ecm{\end{mat}}
\newcommand{\vu}{\vec{u}}

\newcommand{\vpsi}{\vec{\psi}}
\newcommand{\vPsi}{\vec{\Psi}}
\newcommand{\grad}{\nabla}

\newcommand{\T}{{\cal T}}
\newcommand{\hatT}{{\cal \hat T}}
\newcommand{\rhs}{b}  
\newcommand{\dx}{\partial_x}
\newcommand{\dy}{\partial_y}
\newcommand{\dr}{\partial_r}

\numberwithin{equation}{section}

\title{Implicit Adaptive Mesh Refinement for Dispersive \\ Tsunami
Propagation\footnote{Version of March 19, 2024. Revised and resubmitted for publication.}}

\author{Marsha J. Berger\thanks{
Flatiron Institute, 162 5th Ave., NY, NY 10010,
\texttt{mberger@flatironinstitute.org}. Also 
Courant Institute, New York University, 251 Mercer St., NY,
NY 10012, U.S.A. }
\and
Randall J. LeVeque\thanks{Department of Applied Math, University of Washington,
Seattle, WA 98195, U.S.A. Also HyperNumerics LLC, Seattle, Washington,
\texttt{rjl@uw.edu}.}}

\date{December, 2022}
\begin{document}

\maketitle

\begin{abstract}
We present an algorithm to solve the dispersive depth-averaged
Serre-Green-Naghdi (SGN) equations using patch-based adaptive mesh refinement.
These equations require adding additional higher derivative terms to the
nonlinear shallow water equations.
This has been implemented as a new component of the open source GeoClaw
software that is widely used for modeling tsunamis, storm surge,
and related hazards, improving its accuracy on shorter wavelength phenomena.
\add{We use a formulation that requires solving an elliptic system of
equations at each time step, making the method implicit.}
\remove{The equations require the solution of an elliptic system at each
time step.} The adaptive algorithm allows
different time steps on different refinement levels, and solves the implicit
equations level by level. Computational examples are presented to illustrate the
stability and accuracy on a radially symmetric test case and two realistic
tsunami modeling problems, including a hypothetical asteroid impact creating
a short wavelength tsunami for which dispersive terms are necessary.

\end{abstract}



\section{Introduction}\label{sec:intro}

Tsunami propagation and inundation is frequently modeled using the depth-averaged
nonlinear shallow water equations (SWE).  The use of 
adaptive mesh refinement (AMR) is often crucial for realistic problems,
since cell sizes of several km can often be used in the
deep ocean, while inundation of coastal regions generally
requires a 
horizontal resolution of 10 meters or less.  
\add{In our version of patch-based mesh refinement, a nested set of 
finer patches are created and superimposed on the original coarse domain,
time-stepping occurs with an appropriate time step for each level
patch, and regridding occurs every few time steps so that the
finer level patches can track the phenomena of interest
in a wave-propagation problem.
Several figures in this paper show outlines of the refined patches, 
providing a better idea of the how the algorithm works.}
\add{The use of patch-based AMR in computational fluid dynamics has been
well-documented in \cite{Berger-Colella:2-D_AMR, mjb-rjl:amrclaw,
meakin1995efficient,
hornung2002managing,mjb-rjl:amrclaw}; see also
\cite{plewa2005adaptive} for an extensive list of software
packages and applications.  
Our implementation for the shallow water
equations is called GeoClaw and is distributed as
part of the open source Clawpack software \cite{clawpack,MandliEtAl:PeerJ2016}.
It uses high-resolution Godunov-type explicit finite volume methods with
AMR, as described in detail in \cite{BergerGeorgeEtAl2011a,GEORGE20083089,
LeVequeGeorgeEtAl2011}. This software is widely used for tsunami
modeling (among other applications; see \cite{geoclaw-webpage})
and has been accepted as a validated model by the U.S.\
National Tsunami Hazard Mitigation Program (NTHMP) after conducting multiple
benchmark tests as part of
an NTHMP benchmarking workshop \cite{nthmp-benchmarks:2011}.}

The SWE are a hyperbolic system of equations, for which explicit methods work
very well.   These equations are a {\em long wavelength}
approximation, where the underlying assumption is that the water depth
is much smaller than the horizontal length scale.  This is generally suitable
for large-scale earthquake-generated tsunamis, for which the wavelength can be
50 to 100s of km, whereas the typical ocean depth is 4 km.
But for some earthquake generated tsunamis dispersive
equations are more appropriate, and this is almost always true when modeling
phenomena with smaller length scales,
such as tsunamis generated by landslides or asteroid impacts. 
\add{This has been discussed in great detail in previous papers on this
topic, several of which show comparisons between simulations with SWE
and dispersive models, e.g. 
\cite{BabaTakahashiEtAl2015,NASA-AGT-2018,Glimsdal2013,Kim2017,
kirby2022validation, kirby_dispersive_2013, lovholtLynettPedersen_2013,
lovholt2015characteristics, LovholtGlimsdal2015, MatsuyamaIkenoEtAl2007,
Popinet:2015,schambach2019landslide,WardAsphaug:2000,zhou2011nested}.
}
A depth-averaged system of equations is still important for
efficient trans-oceanic propagation. But instead of the shallow water
equations we turn to a Boussinesq-type equation, which retains more terms
from depth-averaging the three-dimensional fluid dynamics. 
These equations, \add{described in more detail below,}
involve second and third-order derivatives, and stability constraints typically
require solving an implicit system of equations each time step.
\add{As a result, solving them can be much more
expensive than solving SWE, for which explicit time-stepping can be
efficiently used. As with SWE, 
the use of AMR can dramatically decrease the computational
expense, and becomes even more important when solving problems with
short-wavelength waves, since much higher resolution may be required in some
regions to capture these waves accurately.
But AMR also becomes much more complex to implement when implicit
systems must be solved on the hierarchy of grids that do not cover the full
domain at most refinement levels.}

The goal of this work is to extend patch-based mesh refinement,
\add{as implemented in GeoClaw,}
to incorporate the solution of an elliptic system each time step.
\remove{This is done in the context of the open source GeoClaw software
\cite{BergerGeorgeEtAl2011a,LeVequeGeorgeEtAl2011} that is
widely used for modeling tsunamis.  It is distributed as part of Clawpack
\cite{clawpack,MandliEtAl:PeerJ2016}
and uses high-resolution Godunov-type explicit finite volume methods with
AMR} \add{By using highly efficient algebraic multigrid methods and
the MPI implementation of PETSc, we have achieved an implementation that can
still solve large-scale realistic tsunami modeling problems on a laptop.}

We implement a form of the 
Serre-Green-Naghdi (SGN) equations developed by Bonneton et al.\
\cite{ti-bo-ma-ch-la:2012},
which is presented in the next section. In
a previous paper \cite{ICM22}, we developed a similar implementation for
another Boussinesq-type system, the Madsen-Sorensen
\cite{MadsenSorensen1992} equations (denoted MS below), but the SGN
equations give a more stable and robust computational method with similar
results.
The elliptic system is formulated and solved separately  on each refinement level,
but includes all patches at that level.
We do this by including the solution variables of the elliptic equation
in the vector of state variables, so that they can
be used as patch boundary conditions (ghost cells) on finer level patches
in the same way as the conserved variables. 
This allows the AMR procedure to incorporate subcycling in time, so that
smaller time steps can be taken on the finer grid patches.
We compare results to those obtained with a composite solver that
does not involve subcycling, in order to verify the accuracy of our
procedure.

There are other Boussinesq solvers that are adaptive in space, including the work
by Popinet \cite{Popinet:2015},
implemented in the Basilisk software \cite{basilisk}.
Our work follows his approach, but with the
added component of refinement in time. 
Other open source tsunami modeling codes that implement dispersive terms
(e.g. \cite{adcirc,coulwave,funwave-tvd,globouss})
do not support general AMR, although many are based on unstructured grids so
that much finer grids are used near coastal regions of interest, or they allow
static nested levels of refinement.
There is also previous work on incompressible flow that uses
adaptive mesh refinement and includes subcycling in time.  The first such work
that we are aware of is the method of Almgren et
al.~\cite{AlmgrenEtAl:1998} for variable density incompressible flow.
Their algorithm used 
an explicit upwind method for the convective terms, 
with the viscous term handled implicitly, and the incompressibility constraint
imposed by solving an elliptic system.
Subsequent recent work by Zeng et al.\ \cite{ZengEtAl:2023} is of the same nature.
Our problem is easier, since our method is less tightly coupled and
does not involve an incompressibility constraint, but it still requires
solving an elliptic system each timestep.  We propose a different algorithm
that also allows for refinement in time.

Near shore it is necessary to switch from the dispersive SGN equations back
to SWE in order to better model wave breaking in very shallow water, and to
robustly handle wetting and drying during coastal inundation. A variety
of criteria have been explored in the literature for optimizing this
transition to best capture wave breaking, \add{e.g.
\cite{kazolea2018wave,KorycanskyLynett,KorycanskyLynett2007,
ti-bo-ma-ch-la:2012}.}
Here we follow \cite{jhkim:phd,Kim2017} and adopt the simplest
approach that is widely used, consisting of using SWE in any grid cell where
the initial water depth in the ocean at rest for that cell or any of its
nearest neighbors is less than some specified tolerance, typically 5
or 10 meters depth.

The grid adaptation and placement
is handled automatically in GeoClaw using a variety of options, such as wave
height, adjoint-based error estimation, or simply
forcing a high level of refinement around particular coastal regions of
interest. The GeoClaw
software incorporates OpenMP for the hyperbolic step by
parallelizing over grid patches.
The implicit system arising in the  Boussinesq equations is solved using an
algebraic multigrid preconditioned Krylov solver in PETSc \cite{balay1998petsc}.
We use the recently introduced PETSc version 3.20 
that allows the use of  MPI for the linear solver in combination 
with the OpenMP parallelization used by GeoClaw. 

This paper is organized as follows.  \Cref{sec:eqn} introduces
the  depth-averaged
SGN equations that have been incorporated into GeoClaw.  We also include the
radially symmetric one-dimensional version of these equations used 
to compute a fine-grid reference solution for some of our 2D test problems.
\Cref{sec:alg} contains a description of the single-grid SGN solver,
and then the block-structured AMR strategy to 
solve the implicit system of equations required by
the SGN solver.  \Cref{sec:ex} contains three computational examples
chosen to validate the stability and accuracy of the implementation and to
illustrate its effectiveness for realistic problems.
We discuss the run times of the new dispersive AMR code and compare them
to similarly refined computations using only the shallow water equations.

\section{Equations}\label{sec:eqn}

The two-dimensional nonlinear shallow water equations (SWE) can be written as
\begin{equation}\label{swe2}
\begin{split}
h_t +(hu)_x +(hv)_y &= 0,\\
(hu)_t + \left(hu^2\right)_x + (huv)_y + gh\eta_x &= 0,\\
(hv)_t + (huv)_x + \left(hv^2\right)_y + gh\eta_y &= 0,\\
\end{split}
\end{equation}
where the surface elevation $\eta(x,y,t) = h(x,y,t) + B(x,y)$, 
with $h$ the fluid depth and $B(x,y)$ the bottom topography.
These are the equations solved in GeoClaw in a manner that handles the
nonlinearity of wave breaking and on-shore inundation very robustly.
However, these equations
are non-dispersive; the linearized equations have the dispersion relation
$\omega(k) = k\sqrt{g h_0}$, and hence constant wave speed
$\omega(k)/k = \sqrt{g h_0}$  for all wave
numbers $k$.

\add{The equations \eqn{swe2} are written in Cartesian (planar) coordinates
and additional terms are added for solving real-world problems on the sphere. 
Coriolis terms can also be added on the sphere, and GeoClaw supports this
option, but many experiments both with GeoClaw and by other researchers have
concluded that Coriolis terms have little effect on global tsunami
propagation due to the very small fluid velocities in the deep ocean 
\cite{Glimsdal2013,kirby2009basin}.
Friction terms
are also added to the momentum equations to model bottom friction based on a
tunable Manning coefficient, as described in
\cite{LeVequeGeorgeEtAl2011} and used by many
standard tsunami modeling codes.  
We omit all these complications in this paper since they are all handled in
the same manner as for SWE in the GeoClaw code.
} 

A number of depth-averaged Boussinesq-type equations have been developed
over the past several decades to model dispersive wave propagation, see
\cite{ma-sc:1998} for one review.
\add{These are generally derived by keeping the next order terms in
asymptotic expansions
in powers of the ratio of water depth to wavelength.}
Different variants are often compared by computing the 
dispersion relation for the linearized version of the equations
and seeing how well it models that of the Airy solution, which is derived
based on the assumption of an incompressible and irrotational flow,
allowing the full three-dimensional velocity to vary with depth.  
The dispersion relation for the Airy solution is
\begin{equation}\label{omegaAiry}
\omega(k) = k\sqrt{gh_0\tanh(kh_0)/kh_0}
\end{equation}
where $k$ is the wavenumber for a wave with spatial wavelength $2\pi/k$ and
$\omega$ is the temporal frequency.

In the previous work of Kim et al.\ \cite{jhkim:phd,Kim2017},
GeoClaw was extended to solve
the MS equations introduced by Madsen and Sorenson \cite{MadsenSorensen1992}
by introducing an implicit solver for an elliptic equation
at each time step. 
The resulting code, BoussClaw, has been used for
solving several dispersive landslide-generated
tsunami modeling problems (e.g.
\cite{gylfadottir_2014_2017,jhkim:phd,lovholt_giant_2017,kim_landslide_2019,liu_modelling_2021}).
However, that extension was implemented only for a single grid patch at a
fixed resolution, and did not use the AMR aspects of GeoClaw.

In our work on extending the AMR algorithms to work for Boussinesq
equations, we started with the BoussClaw code and our first AMR implementation
again solved the MS equations.  The algorithms and some sample results
were presented in \cite{ICM22}.  However, these algorithms exhibited
poor stability properties at patch interfaces and non-reflecting boundaries,
and numerous attempts to make
them more stable have not resulted in a sufficiently robust code.  
We now believe there are inherent stability issues with these equations,
which will be explored further elsewhere.
L{\o}vholt and Pedersen \cite{LovholtPedersen2009} also observed instabilities
for similar equations when modeling flow over topography. 

We have switched to using a different form of Boussinesq equations,
an improved form of the Serre-Green-Naghdi (SGN) equations developed in
\cite{ti-bo-ma-ch-la:2012}, which includes a parameter $\alpha$ that can be tuned
to match the Airy dispersion relation quite well.
The dispersion relation for this improved SGN equation is
\begin{equation}\label{omegaSGN}
\omega(k) = k\sqrt{\frac{gh_0(1 + (\alpha-1)(kh_0)^2/3)}{1 + \alpha(kh_0)^2/3}},
\end{equation} 
Using $\alpha=1$ gives the original SGN equations, but 
following the choice made by Popinet in the Basilisk code we use
$\alpha = 1.153$ (which gives an imperceptible change in practice from the 
value 1.159 originally suggested in \cite{ti-bo-ma-ch-la:2012}).

\begin{figure}
\hfil\includegraphics[width=5.5in]{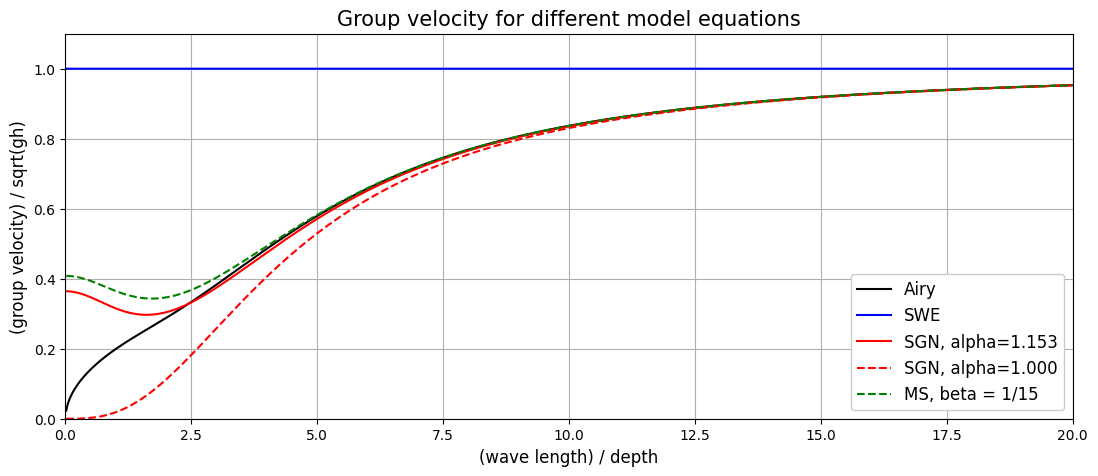}\hfil
\caption{\label{fig:gpvel} 
The scaled group velocity $\omega'(k)/\sqrt{gh_0}$ for the dispersion
relations shown in the text, plotted versus $2\pi/(kh_0)$,
the wavelength relative to the undisturbed fluid depth.
In this work we implement SGN with $\alpha=1.153$.
  }
\end{figure}

For comparison, we also note that the MS equations have a tunable
parameter $\beta$ and dispersion relation
\begin{equation}\label{omegaMadsen}
\omega(k) = k\sqrt{\frac{gh_0(1 + \beta(kh_0)^2)}{1 + (\beta+1/3)(kh_0)^2}},
\end{equation} 
with the value $\beta = 1/15$ giving the best agreement with the Airy
solution. Setting $\beta = 0$ in \eqn{omegaMadsen} or $\alpha=1$ in
\eqn{omegaSGN} gives the dispersion relation of the original equations of
Serre-Green-Naghdi.  \Cref{fig:gpvel} shows the group velocity
$\omega'(k)$ plotted versus $2\pi/(kh_0)$ (the wavelength divided by the
fluid depth) for each set of equations, scaled by $\sqrt{gh_0}$.
Note that for wavelengths greater than approximately twice the fluid depth,
the SGN equations with $\alpha=1.153$ give a very good approximation to the
Airy solution in the linearized case.  
In the long-wave limit all of the models asymptote
towards the linearized shallow water limit.  For very short wavelengths
the Airy group velocity approaches 0.
So does the group velocity for SGN with $\alpha=1$, but much too slowly, often
resulting in highly oscillatory waves that trail far behind the main wave in
an unrealistic manner.  By contrast, SGN with the value of $\alpha = 1.153$
has group velocities bounded away from zero. This is
evident in some of the examples of \cref{sec:ex}, where radially expanding
waves leave a quiescent central zone in their wake, which is generally
preferable to highly oscillatory waves that are not
physically meaningful.

The SGN system of equations has the form of the SWE \eqn{swe2} with
the addition of source terms:
\begin{equation}\label{sgn2}
\begin{split}
h_t +(hu)_x +(hv)_y &= 0\\
(hu)_t + \left(hu^2\right)_x + (huv)_y + gh\eta_x &= 
    h\left(\frac{g}{\alpha} \eta_x - \psi_1\right) \\
(hv)_t + (huv)_x + \left(hv^2\right)_y + gh\eta_y &= 
    h\left(\frac{g}{\alpha} \eta_y - \psi_2\right)
\end{split}
\end{equation}
The vector $\vpsi = [\psi_1,~\psi_2]$ satisfies an elliptic equation of the
form
\begin{equation}\label{sgn1}
(I + \alpha {\cal T})\vpsi = \rhs
\end{equation} 
where $\vu = [u,v]$ is the velocity vector.  
These equations agree with equation (4) of \cite{Popinet:2015}, but with a typo
corrected (there is no $h$ on the left hand side).
The second-order elliptic operator ${\cal T}$ and right hand side $\rhs$ are
written in vector notation in \cite{Popinet:2015}. Here we write them out in
component form to clarify the derivatives that must be approximated in the
finite difference implementation of these equations.

The components of the differential operator
\[
\T = \bcm \T_{11}&\T_{12} \\ \T_{21}&\T_{22} \ecm
\]
are given by
\begin{equation}\label{T}
\begin{split}
\T_{11} &= -\frac{h^2}{3} \dx^2 - hh_x \dx
    + \left(\frac h 2 B_{xx} + B_x\eta_x\right), \\
\T_{12} &= -\frac{h^2}{3} \dx\dy + \frac{h}{2}B_y\dx 
    -h\left(h_x + \frac 1 2 B_x\right) \dy
    + \left(\frac{h}{2}B_{xy} + B_y\eta_x\right), \\
\T_{21} &= -\frac{h^2}{3} \dx\dy 
    -h\left(h_y + \frac 1 2 B_y\right) \dx
    + \frac{h}{2}B_x\dy 
    + \left(\frac{h}{2}B_{xy} + B_x\eta_y\right), \\
\T_{22} &= -\frac{h^2}{3} \dy^2 - hh_y \dy
    + \left(\frac h 2 B_{yy} + B_y\eta_y\right). \\
\end{split}
\end{equation} 
The right hand side $\rhs$ in \eqref{sgn1} has components
\begin{equation}\label{rhs}
\begin{split}
\rhs_1 &= \frac{g}{\alpha} \eta_x +2h\left(\frac h 3 \phi_x 
    + \phi \left(h_x + \half B_x\right)\right)
    + \frac h 2 w_x + w\eta_x, \\
\rhs_2 &= \frac{g}{\alpha} \eta_y +2h\left(\frac h 3 \phi_y
    + \phi \left(h_y + \half B_y\right)\right)
    + \frac h 2 w_y + w\eta_y,
\end{split}
\end{equation} 
where
\begin{equation}\label{phiw}
\begin{split}
\phi &= \dx\vu\cdot\dy\vu^{\perp} + (\grad\cdot\vu)^2\\
     &= v_xu_y - u_xv_y + (u_x+v_y)^2,\\
\noalign{\vskip 5pt}
w &= \vu\cdot(\vu\cdot\grad)\grad B\\
  &= u^2B_{xx} + 2uvB_{xy} + v^2B_{yy}.
\end{split}
\end{equation}

\subsection{Radial symmetry}\label{sec:radial}
For testing purposes it is very convenient to consider problems in which the
initial data and bathymetry are radially symmetric, in which case the
Boussinesq equations reduce to equations in one space dimension $r$ and time.
Numerical solutions obtained with the
two-dimensional code can then be compared with fine-grid solutions of the radial
equations to verify that the complicated matrix equations have been properly
implemented.  This is done in \cref{sec:radocean} below.

We use $\hat h(r,t)$ to denote the depth in radial coordinates and similarly
for other variables.  The radial velocity $\hat u(r,t)$ corresponds to the
two-dimensional velocity field $u(x,y,t) = \hat u(r,t)\cos(\theta), 
v(x,y,t) = \hat v(r,t)\sin(\theta)$.  Converting the Boussinesq system to
polar coordinates and assuming there is no variation in the $\theta$
direction, we obtain the equations

\begin{equation}\label{sgnradial}
\begin{split} 
\hat h_t + \frac 1 r \left(r\hat h\hat u\right)_r &= 0,\\
(\hat h \hat u)_t + \frac 1 r \left(r \hat h\hat u^2\right)_r 
+ g\hat h \hat\eta_r &= \hat h \left(\frac g \alpha \hat\eta_r -
\hat\psi\right),
\end{split} 
\end{equation} 
where $\hat\psi$ satisfies the one-dimensional elliptic equation
\begin{equation}\label{ell1}
(1 + \alpha \hatT)\hat\psi = \hat b.
\end{equation} 
In this case the second-order scalar differential operater $\hatT$ is given by
\begin{equation}\label{hatT}
\hatT = -\frac{\hat h^2}{3} \left(\dr^2 + \frac 1 r \dr - \frac{1}{r^2}\right)
- \hat h\hat h_r\left(\dr + \frac 1 r \right) 
+ \frac 1 2 \hat h \left(\hat B_{rr} - \frac 1 r \hat B_r\right) + \hat B_r\hat\eta_r.
\end{equation}
The right hand side  of \eqn{ell1} is
\begin{equation}\label{hatb}
\hat b = \frac{g}{\alpha}\hat\eta_r + 2\hat h\left(\frac{\hat h}{3}
\hat\phi_r + \hat\phi\left(\hat h_r + \frac 1 2 \hat B_r\right)\right)
+ \frac{\hat h}{2}\hat w_r + \hat w\hat\eta_r,
\end{equation} 
where, analogous to \eqn{phiw}, $\hat\phi$ and $\hat w$ are given by
\begin{equation}\label{hatphiw}
\begin{split} 
\hat\phi &= (\hat u_r)^2 + \frac 1 r \hat u_r\hat u + \frac{1}{r^2}\hat u^2,\\
\hat w &= \hat u^2 \hat B_{rr}.
\end{split}
\end{equation} 

Note that far from the origin the radial waves should behave like plane
waves, and dropping the terms that depend on $1/r$ from these equations
reduces them to the plane wave case.  This consists of the
one-dimensional shallow water equations with a source
term $\psi_1$ in the momentum equation that comes from solving the
one-dimensional elliptic equation $(I+\alpha\T_{11})\psi_1 = b_1$, where
$\T_{11}$ and $b_1$ are given by \eqn{T} and \eqn{rhs}, respectively.

\section{Solution Algorithm}\label{sec:alg}
We give a high level overview of our approach before giving the details
below.
Following Popinet and others, we use a splitting method to take a full
time step. The elliptic equation \eqn{sgn1} is discretized using a 
straightforward second order
finite difference scheme, and the nonlinear shallow water equations are
solved in finite volume form as implemented in GeoClaw.
A time step starts by solving the elliptic equation, which updates
only the momenta.
Then the shallow water equations are solved for one time step using 
the modified momenta and the original depth $h$ 
to get to the next time level.
The details of this splitting method are given in \cref{sec:alg1}.
In \cref{sec:alg2} 
we describe the patch-based AMR version of the splitting method. 
This is the major algorithmic innovation in this paper. 
We solve an  elliptic equation one level at a time, but
including all patches at the given level. We also introduce a second
``provisional'' solve on coarser levels at the next time step. 
This allows linear interpolation
in time to fill in the ghost cells on finer levels, which is needed for 
subcycling in time.
Incorporate the elliptic solution variables from each level into the state
vector made this easy to implement in the existing AMR framework.
In \cref{sec:alg3} we describe a composite solve of a linear system
coupling all levels, and without
subcycling in time, which we use for comparison.

We use capital letters to denote the numerical solution on a single grid
or the coarse grid; below we introduce lower case for the fine grid solution.
To enable the subcycling in the AMR implementation, the vector of variables 
$\vPsi = (\Psi_1, \Psi_2)$
from the implicit solve are stored along with the conserved variables
$(H,HU,HV)$, defining $Q = (H,HU,HV,\Psi_1,\Psi_2)$. 
The elliptic equations for the dispersive update are solved over the entire
domain except when the water depth falls below a threshold. This allows
the SWE to handle the wetting and drying on land. They are also better
at simulating wavebreaking, which does not naturally occur in the
dispersive equations.  For each of the examples in 
\cref{sec:ex} we specify the undisturbed water depth
used as the threshold for switching from SGN to SWE. 

\subsection{Splitting Method for SGN}\label{sec:alg1}

A single time step of the fractional step algorithm on a single grid proceeds
as follows. We start with $(H,HU,HV)^N$ at time $t^N$. On a single
grid there is no need to store the $\vPsi^N$ vector, so we simplify the
description here by omitting it.  The following
steps advance the solution to $t^{N+1} = t^N + \Delta t$:
\begin{enumerate}
\item
Solve the elliptic equations for $\vPsi^N$ in \eqref{sgn1}. 
The right-hand side and
matrix coefficients are a function of $(H,HU,HV)$ 
at time $t^N$. 

We currently use PETSc to solve the sparse system of equations. The
matrix elements (in compressed sparse row form) along with the right-hand
side vector are passed to the algebraic multigrid preconditioned 
Krylov solver in PETSc.
\remove{ with a change in the default convergence criteria
from $10^{-5}$ to $10^{-7}$. 
We reuse the preconditioner for multiple time steps, which
approximately halves the run time. }

Popinet \cite{Popinet:2015} and others have implemented a multigrid
solver and report good convergence results. Since they use a fixed spatial
refinement of a factor of 2, the multigrid hierarchy fits nicely into their
strategy. We often refine by larger factors, and have not gone
this route.

\item
Advance the momentum equations using forward Euler and the
right-hand-side of \eqref{swe2} to get 
\begin{equation}
\begin{split} 
H^* &= H^N,\\
(HU)^* &= (HU)^N + \Delta t \, H^N\left(\frac{g}{\alpha} \eta^N_x - \Psi^N_1\right),
 \\
(HV)^* &= (HV)^N + \Delta t  \, H^N\left(\frac{g}{\alpha} \eta^N_y - \Psi^N_2\right). 
\end{split} 
\end{equation}
In these equations $\eta_x$ and $\eta_y$ now stand for centered finite
difference approximations to these quantities.
We have also experimented with a 2-stage Runge-Kutta method but found
almost no difference, so it was not worth the computational expense.
\add{This is likely due to the first order errors that are introduced
with limiters and the splitting error between the shallow
water and dispersive steps. Moreover the SWE already have reduced
accuracy with realistic bathymetry.}
\item
Take a time step $\Delta t$ with the SWE solver, with initial data
$(H,HU,HV)^*$ to obtain the values $(H,HU,HV)^{N+1}$.
We denote this update by $(H,HU,HV)^{N+1} = SW((H,HU,HV)^*, \dt)$.
\end{enumerate}

At domain boundaries, the conditions for $\vPsi$, e.g.
extrapolation or wall boundaries, are put
directly into the matrix equations. For example, the equation for
a cell $(i,j)$ adjacent to the right edge of the computational domain would
normally include cell $(i+1,j)$ in its stencil.
To implement wall boundary conditions, we want to specify
$\Psi_{1,i+1,j} =  -\Psi_{1,i,j}$ since this is a correction to the
momentum in the $x$-direction.
We implement this by negating the matrix coefficient that would correspond 
to the missing cell, and adding it to cell $(i,j)$'s matrix entry. 
Since $\Psi_2$ updates the $y$-momentum, tangential to this edge, we want
$\Psi_{2,i+1,j} =  \Psi_{2,i,j}$ and the matrix coefficient is not
negated before adding it to cell $(i,j)$'s matrix entry.
At the top and bottom boundaries, wall boundary conditions are
implemented similarly with the negation swapped between ${\Psi_1}$ and
${\Psi_2}$.
For extrapolation boundary conditions, no components are negated.
\add{Extrapolation boundary conditions do not absorb outgoing waves as
well for SGN as for SWE, where they are routinely used in GeoClaw at
open ocean boundaries.  The implementation of better absorbing
boundary conditions is still work in progress, using a sponge layer
as often used in for other dispersive solvers, e.g. \cite{LovholtGlimsdal2015}.
Boundary conditions for dispersive equations are still an open 
area of research \cite{noelleEtAl_bcs}.}

The switch from SGN to SWE is made on a cell-by-cell basis.
If the inital water depth for any cell
in the 3-by-3 neighborhood surrounding a given cell is below a specified
threshold, then the dispersive correction
is not applied.  This is implemented by setting this
row of the linear system for $\vPsi$ to the corresponding
row of the identity matrix and the
right-hand side to zero, ensuring a zero correction.
The cell can still be used in the stencil for the dispersive terms
of a neighboring cell if appropriate.
We have not observed any spurious reflections from the interface where the
equation changes when using this procedure.

\subsection{AMR for SGN}\label{sec:alg2}
\begin{figure}[thb]
\centerline{\includegraphics[height=2.3in]{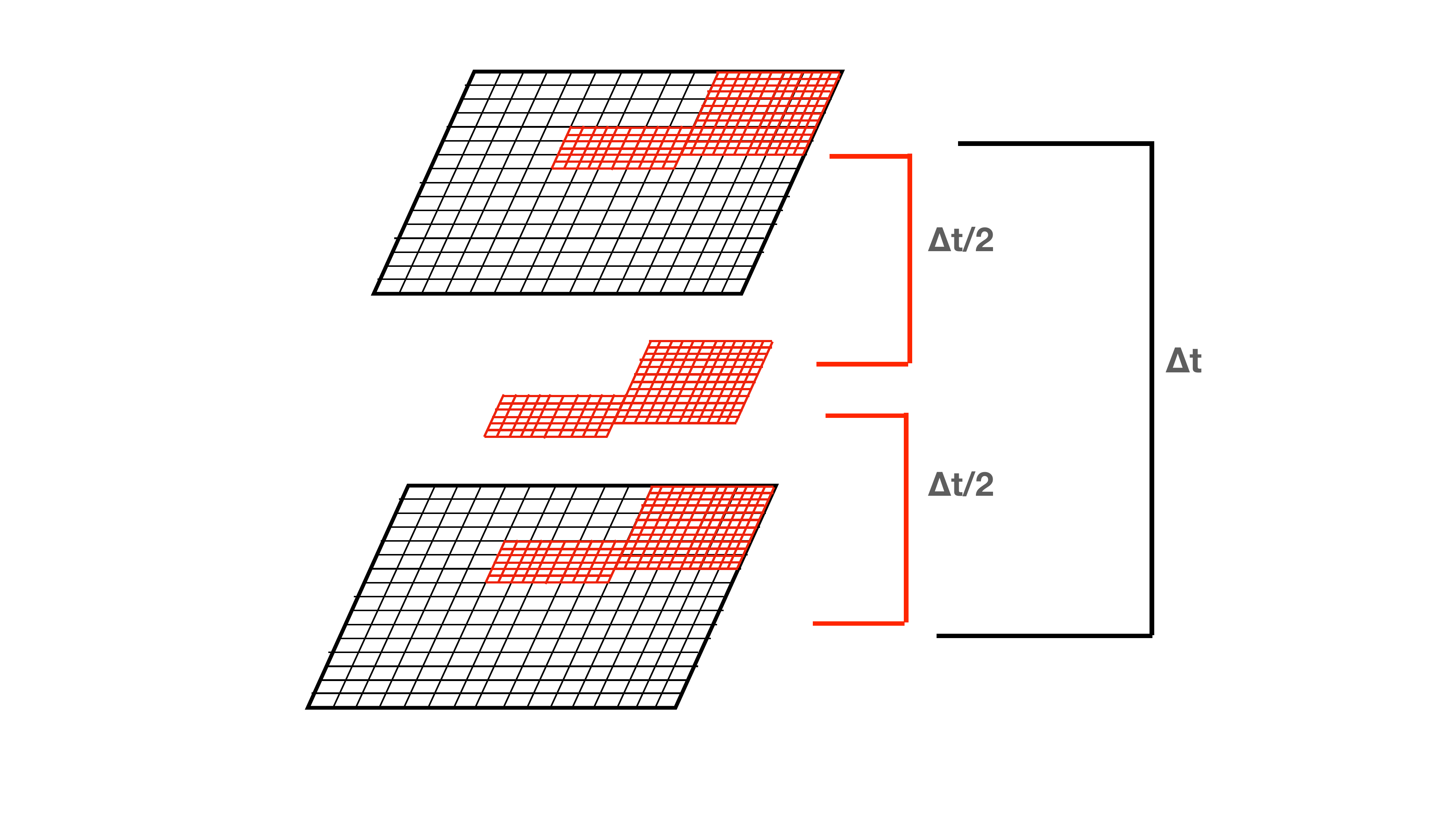}}
\caption{There is a coarse grid with cells outlined in black on the base
grid. Over part of the domain there are
two adjacent fine patches refined by a factor of two, with cells
outlined in red. The fine grid time step is half the coarse grid step. Before
the fine grid takes a step, two ghost cells (not shown in the figure)  are 
needed to complete the stencil. Figure taken from \cite{ICM22}.\label{fig:amrfig}}
\end{figure}

Now suppose we have two grid levels with refinement by a factor of 2 in time.
The main changes to the algorithm involves fine grid ghost cells, and a
second implicit solve on the coarser levels.
\Cref{fig:amrfig} illustrates the setup, with one grid at the
coarsest level, and two adjacent rectangular non-overlapping fine grid patches at level 2.
We denote the coarse grid values at some time $t^N$ as above, using
capital letters. 
We assume that the fine grid is at time $t^N$ , but that on the 
fine grid we want to take two time steps of $\Delta t/2$ to reach time
$t^{N+1}$.
For the fine grid values at time $t^N$ we use
lower case letters with $q^N=(h,hu,hv,\psi_1,\psi_2)^N$.

We always need boundary conditions in ghost cells around the union of
fine grid patches when taking an explicit hyperbolic time step, which in
general are obtained by space-time interpolation from the underlying
coarser level grids.  For the SGN equations we also need boundary conditions
for the elliptic system of equations that is solved on the union of all
grid patches at the fine level. By incorporating $\vPsi$ in $Q$, the same
space-time interpolation operators provide these necessary ghost cell values.

In the algorithm description below, $\If(Q^N)$ denotes the spatial
interpolation operator that
interpolates from coarse grid values to the ghost cells of a fine grid
patch at a single time $t^N$.

One time step on the coarse grid, coupled with two time steps on the
fine grid, is accomplished by the following steps:

\begin{enumerate}

\item {\bf Coarse grid step:}
    \begin{enumerate}
    \item Take time step $\Delta t$ on the coarse grid as described above for
    the single grid algorithm, but denote the result by
    $(\widetilde H,\widetilde{HU},\widetilde{HV})^{N+1}$ since these 
    provisional values will later be updated.
    \add{This step is taken on the union of all coarse level patches,
    and does not involve fine patches.}
    \item  Solve for
    a provisional $\widetilde\Psi^{N+1} = [\Psi_1,\Psi_2]$. This is the
    second implicit solve mentioned above. This will be needed for interpolation
    in time when determining ghost cell values for $\psi$ on the fine grid
    using $\If(Q^N)$ and $\If(\widetilde Q^{N+1})$.
    \end{enumerate}

\item {\bf Fine grid steps:}
    \begin{enumerate}
    \item Given $(h,hu,hv)^N$ and ghost cells boundary 
    conditions $\If(Q^N)$ on the union of fine grid patches, solve the
    fine grid elliptic system for $\psi^N$.
    \item Update using the source terms 
    $$
    (h,hu,hv)^* = (h,hu,hv)^N +
    \left(0, \frac{\dt}{2}\left(\frac{g}{\alpha} \eta_x - \psi_1 \right)^N,
    \frac{\dt}{2}\left(\frac{g}{\alpha} \eta_y - \psi_2 \right )^N\right) .
    $$
    \item Take a shallow water step: 
   $$
   (h,hu,hv)^{N+1/2} = SW((h,hu,hv)^*, \dt/2).
   $$
    Note that we use $t^{N+1/2} = t^N + \dt/2$ to denote the intermediate time
    for the fine grid.
    \item Obtain ghost cell values at this intermediate time as
    $\half(\If(Q^N)+\If(\widetilde Q^{N+1}))$.
    \item Solve the elliptic system for $\psi^{N+1/2}$ on the fine grid.
    \item Update using the source terms 
    $$
    (h,hu,hv)^* = (h,hu,hv)^{N+1/2}  +
    \left(0, \frac{\dt}{2}\left(\frac{g}{\alpha} \eta_x - \psi_1
    \right)^{N+1/2},
    \frac{\dt}{2}\left(\frac{g}{\alpha} \eta_y - \psi_2 \right
    )^{N+1/2}\right).
    $$
    \item Take a shallow water step: 
    $$
       (h,hu,hv)^{N+1} = SW((h,hu,hv)^{*}, \dt/2).
    $$
    \end{enumerate}

\item {\bf Update coarse grid:}
    \begin{enumerate}
    \item Define $(H,HU,HV)^{N+1}$ by the provisional values  
    $(\widetilde H,\widetilde{HU},\widetilde{HV})^{N+1}$ where there is 
    no fine grid covering a grid cell, but replacing
    $(\widetilde H,\widetilde{HU},\widetilde{HV})^{N+1}$ by the 
    conservative average of $(h,hu,hv)^{N+1}$
    over fine grid cells that cover any coarse grid cell.
    \add{In other words, the fine grid values update all underlying coarser cells to get the
    coarser level ready for its next step.}
    \end{enumerate}
    \end{enumerate}

The final step is applied because the fine grid values $(h,hu,hv)^{N+1}$ are
more accurate than the provisional coarse grid values.
Note that it is not necessary to update $\Psi$ to $\Psi^{N+1}$
on the coarse grid because this is computed at the start of the next time step.

We then proceed to the next coarse grid time step.
At the start of this step, the updated $(H,HU,HV)^{N+1}$ is used
to solve for $\Psi^{N+1}$, and the provisional $\widetilde \Psi^{N+1}$ is
discarded.
Hence two elliptic solves are required on the coarse level each time step,
rather than only one as in the single grid algorithm.
Luckily this second solve is not needed on the finest grid level.

The algorithm above easly generalizes to refinement factors larger than 2. 
There will be additional time steps on level 2, and for each
time step the ghost cell BCs will be determined by linear interpolation in
time between $\If(Q^N)$ and $\If(\widetilde Q^{N+1})$.
If there are more than two levels, this same idea is applied at the next
level.
After each time step on level 2, any level 3 grids will be advanced by the
necessary number of time steps to reach the advanced time on level 2.  In this
case we will need two elliptic solves for every time step on level 2,
one at the start of a level 2 time step, and one for the provisional
values after advancing level 2, to provide interpolated ghost cell values
to level 3.

\begin{figure}[htb]
\centerline{\includegraphics[height=2.1in]{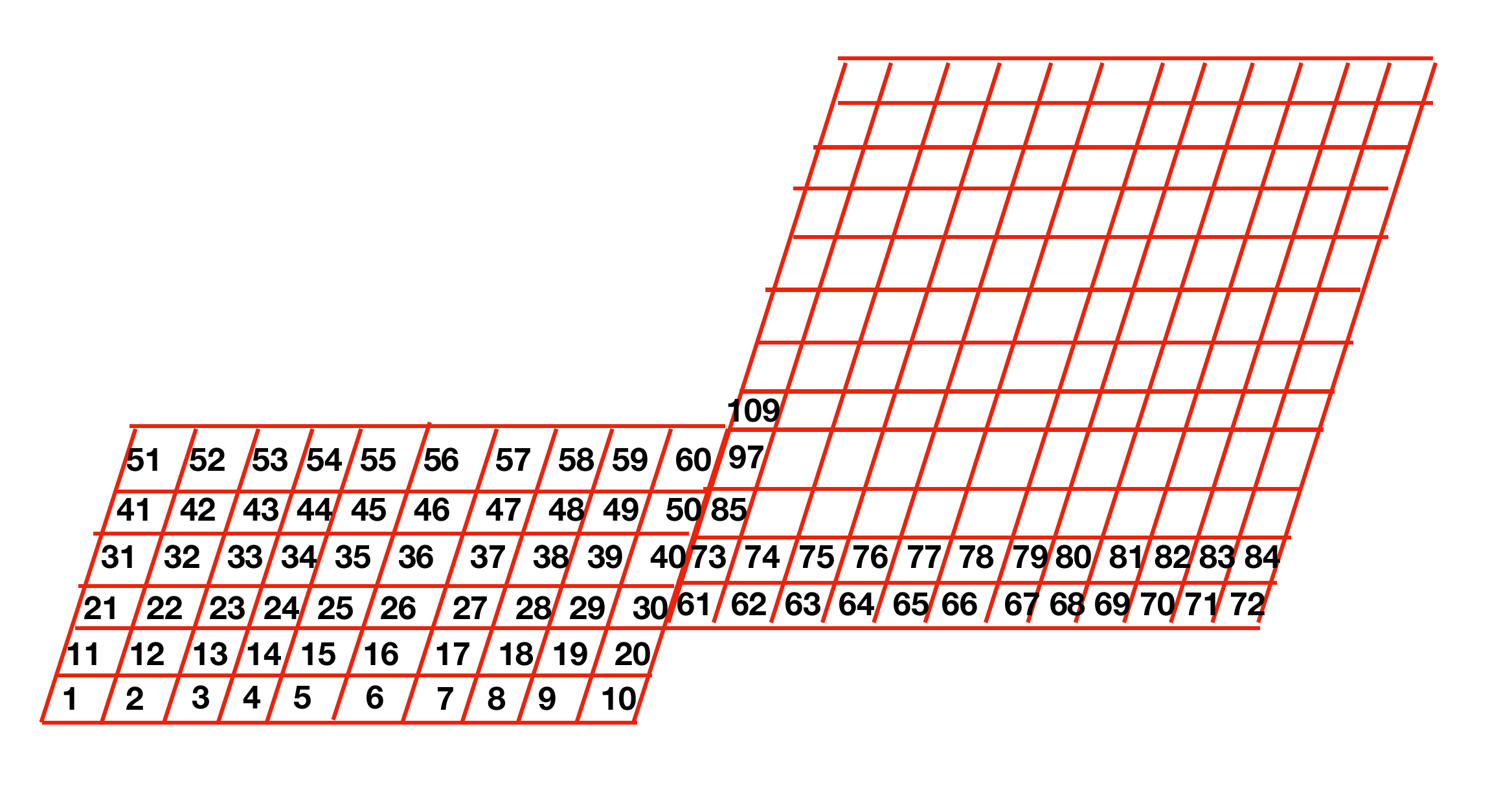}}
\caption{Sample enumeration of cells on  level 2 grids. The
numbers are used to map to the row in the linear system. 
\label{fig:marking}}
\end{figure}

When there are multiple possibly adjacent grids at the same level, as in
\cref{fig:amrfig}, the
border cells of one grid might include another grid's interior cells in its
stencil. Its row and column numbers need to be known to build the entries
of the matrix. To enable this, grids are pre-processed at every adaptive
regridding step by numbering the cells to indicate what equation in the
matrix they apply to. This information is saved on each grid patch. 
\Cref{fig:marking} illustrates this enumeration
for the level 2 grids in \cref{fig:amrfig}.
For example, the stencil for cell 40 on the left grid involves cell 73
from the right grid. The preprocessing step would fill in these numbers in
those ghost cells that come from adjacent grids. 
\remove{These numbers are then
used to build the (row, column, matrix entry) triplets in  the sparse matrix 
one patch at a time.} 
If there are no such adjacent grid patches
then the ghost cells would have been interpolated from the coarse grid,
and are treated as Dirichlet conditions. They would be flagged with $-1$
to indicate that they have no row or column associated with them, and
their values are incorporated into the right-hand side of the equations.

For adaptive calculations, the PETSc preconditioner mentioned above 
can be reused until the grids at that level change. The coarsest level 1
grids never change, but they are so coarse and inexpensive that reuse is
not important.

\add{The finite volume method used for SWE in GeoClaw
exactly conserves mass on a Cartesian
domain with no coastal or onshore points, even when AMR is used.
(And also near the coast on a uniform grid,
but see \cite{LeVequeGeorgeEtAl2011} for
more discussion of conservation and the issues that arise when using
AMR near the coast.)  GeoClaw also conserves momentum when used
on a flat bottom, but with any realistic topography momentum should not be
conserved.  The modifications we introduce to
solve Boussinesq-type equations add additional source terms only 
to the momentum equations, and so mass is conserved as well as it is for SWE.
}

\subsection{Composite Solution Algorithm}\label{sec:alg3}

To verify the accuracy of our procedure for refinement in time, we also 
implemented a multilevel composite solve without subcycling. 
\add{In this composite solve, the unknowns on all levels are coupled
together in a single system of equations. In this approach all levels take
the same time step. }
\remove{All levels
are solved together in one larger coupled system of equations. }
This global time step must then be chosen for stability on the finest
level, based on the CFL condition for the explicit SWE solver.
This typically means that the Courant number is much smaller
than 1 on the coarser levels, requiring more coarse time steps and perhaps
diminishing the accuracy due to numerical dissipation. For this reason we
prefer to use subcycling in practice.

\begin{figure}[htb]
\centerline{\includegraphics[height=2.0in]{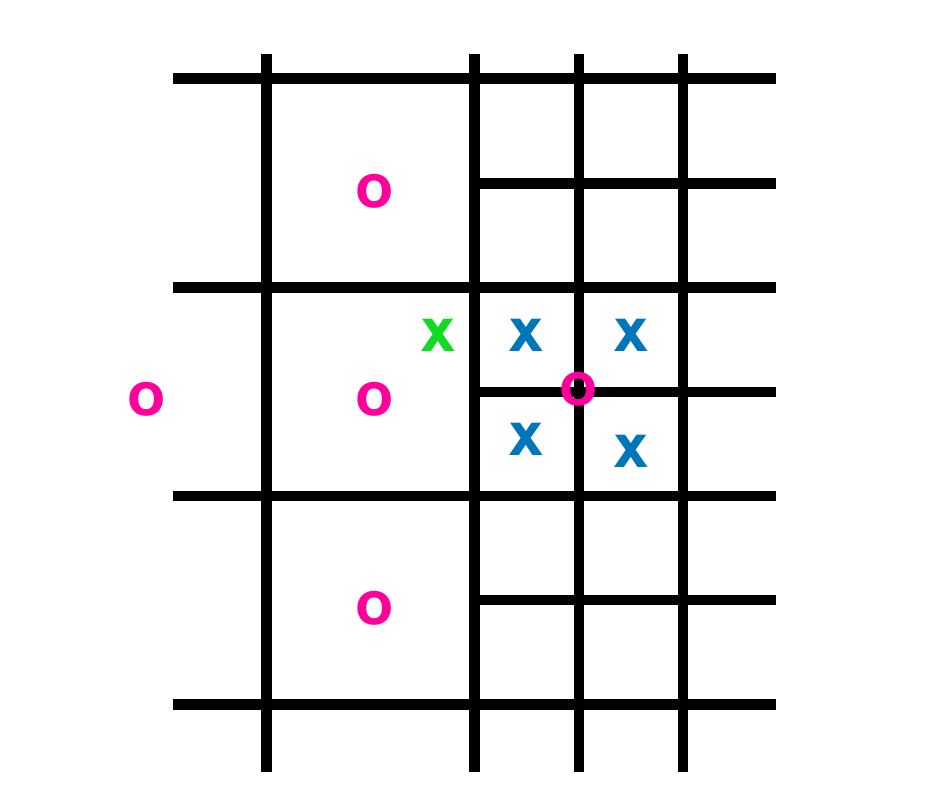}}
\caption{For the composite coupled solution, new unknowns and equations
are introduced into the elliptic solver. The fine
grid ghost cell marked by the green x is a new unknown, and its
equation specifies that this value is linearly interpolated from the 5
coarse grid values indicated by red circles.
Note that the green cell is also part of several fine grid stencils for
the interior fine cells adjacent to the patch edge.
There is also a new unknown on the coarse grid marked with a red circle
that is underneath the fine grid. Its equation specifies
that the coarse grid value is the average of the 4 fine grid values.
\label{fig:intfc}}
\end{figure}

The composite solve introduces two new types of equations, pictorially
illustrated in \cref{fig:intfc}. The first and
last rows and columns just inside the border of
a fine grid have an equation that says their area-weighted sum should equal 
the corresponding value on the underlying coarse cell. 
This equation is associated with the coarse cell index, and it
needs to be included in the grid enumeration procedure.
Note that this {\em hidden} coarse cell is also
part of its coarse neighbor's stencil that lies outside the fine grid.

Secondly, all the ghost cells on the fine grid will now have equations for
$\vpsi$ that   
specify linear interpolation in space only from the coarse grid. 
We use a centrally
differenced gradient in $x$ and $y$ on the coarse cell to reconstruct to
the fine ghost cell. If the
stencil includes a cell that lies outside the computational domain 
the boundary conditions are applied.
Note that the fine grid ghost cell is also referenced in the stencil 
of the first interior fine  cell mentioned above.

Except for the border cells mentioned above, the rest of the coarse cells 
that are
{\em underneath} a fine grid do not need to participate in the composite
solve. They are marked with a flag to ignore them 
in the composite enumeration pre-processing step. After the solve,
the coarse cells are updated by the fine grid, and the shallow water
step proceeds on each level separately as in the usual patch-based AMR
for hyperbolic equations, except with the same time step.

\subsection{Algorithmic comparisons}\label{sec:alg5}
We compare the composite coupled algorithm just described with
the patch-based algorithm with subcycling in time  described in \cref{sec:alg2}.
The problem has a flat topography with water depth 4000 m.  
The radially symmetric initial
sea surface is given by $\eta(x,y,0) = \exp(-(r/2000)^2)$,
where $r = \sqrt{x^2+y^2}$.  The initial velocity is 0 everywhere.
The simulation is run to time 300 seconds.
There are 5 levels of refinement, each
with a factor of 2 in $x$ and $y$. The composite algorithm has a
single $\Delta t$ based on the stability limit for the level 5 grid.
The subcycled algorithm refines by a factor of 2 in time as well.
We also compare with a uniformly fine solution at the same resolution
as the finest level.

\begin{figure}[h!]
\centerline{\includegraphics[width = .9\textwidth]{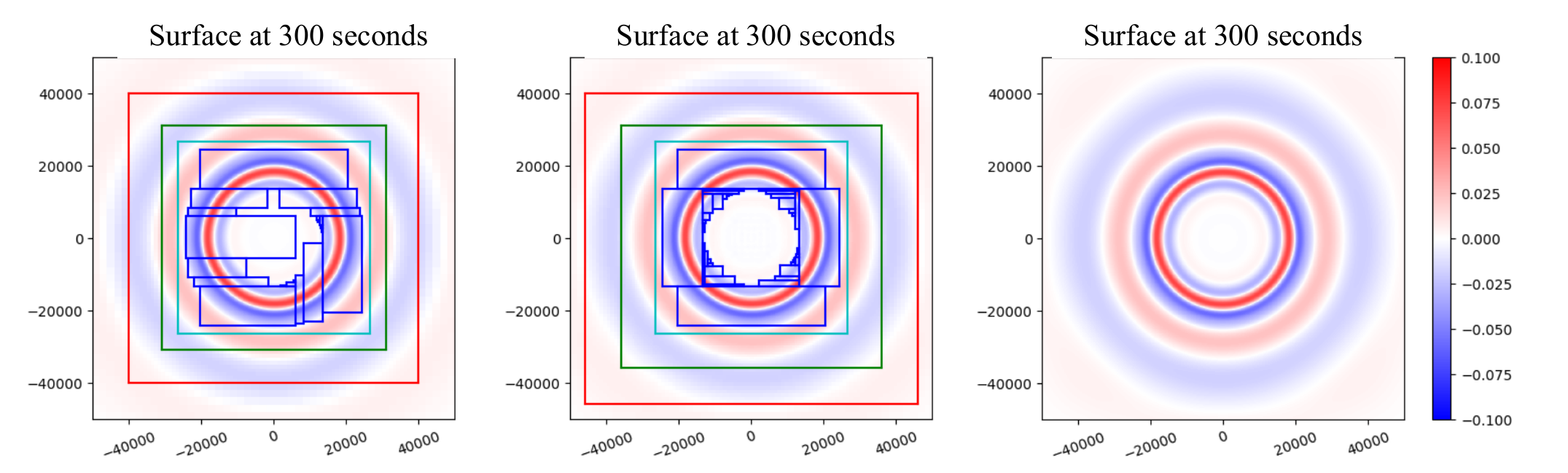}}
\centerline{\includegraphics[width = .7\textwidth]{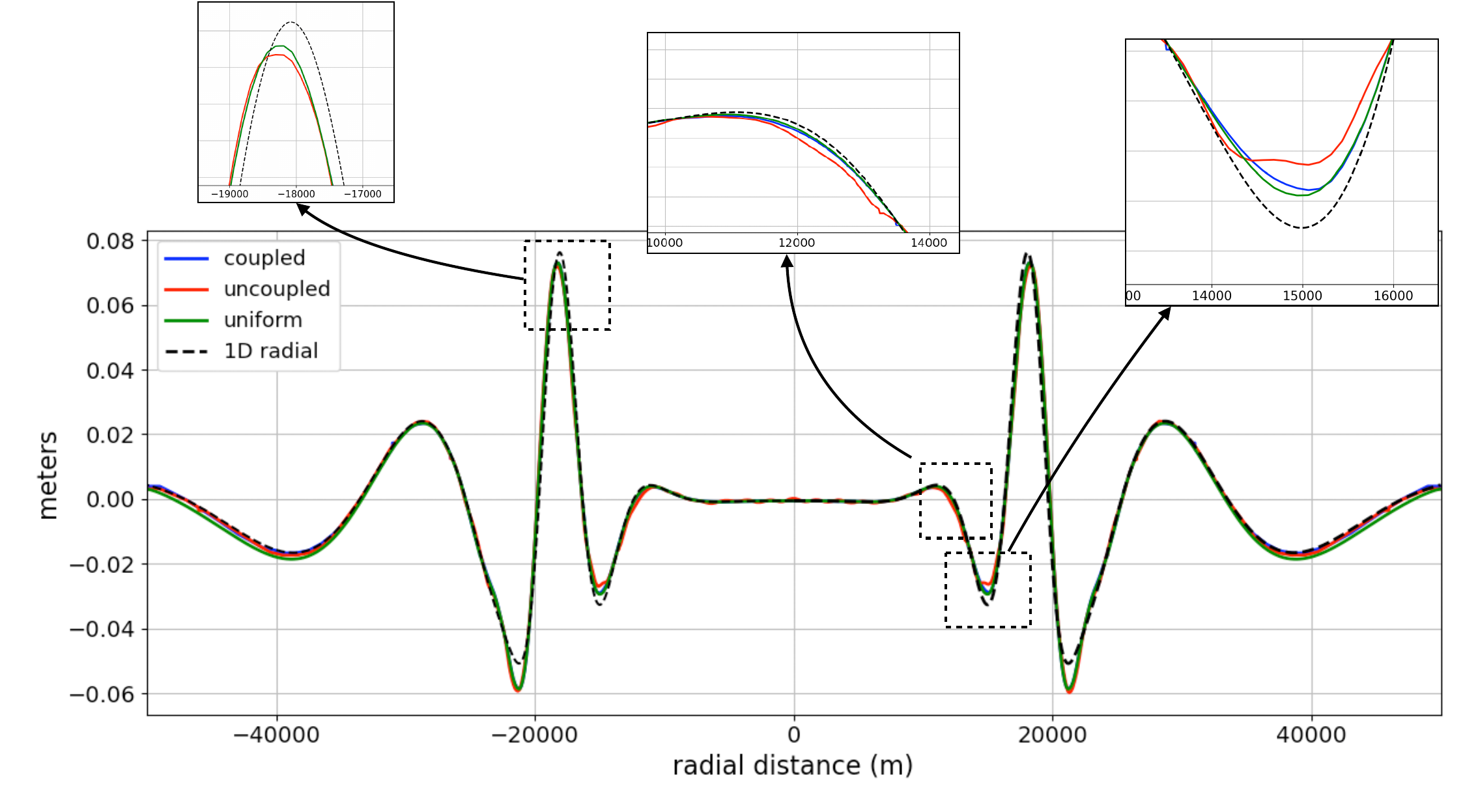}}
\caption{\remove{Comparison of three different runs.}
\add{Comparison of 5-level coupled and uncoupled AMR results with results
on a uniform grid having the same resolution as the finest AMR level.}
Top left is the composite coupled algorithm; middle is the uncoupled algorithm
with subcycling in time; right is the equivalently refined uniform grid.
The grid refinement patterns are fairly different but the
solutions agree very well. Horizontal transects through the  3 solutions are
shown in the lower plot, and compared to the reference 1D radial solution. 
The zoomed sections show
the very small differences in the uncoupled solution (red curve) 
from the others. 
\label{fig:radialEta}}
\end{figure}

\Cref{fig:radialEta} shows the results from the 3 approaches
after 300 seconds.  To the eye the leading waves are identical.
The two adaptively refined
simulations show quite different refinement patterns, none of which enforce
symmetry, yet the solutions remain remarkably symmetric. 

\add{For completeness, we also compared the coupled algorithm with an
uncoupled solution and the same fixed (small) time-step.}\remove{with a patch-based uncoupled solution
but with no refinement in time (not shown here), usingthe same small fixed time
step. }
\remove{Note that this is not the same as the composite solve
method, which solves all levels together. The subcycled algorithm (even
with no refinement in time) solves an implicit system on level 1 alone,
interpolates for level 2 ghost cells and solves on level 2 alone, etc.}
This intermediate algorithm was also more 
accurate than the subcycled algorithm in a few places, at the cost of  
increased run-time due to the much smaller time step needed here,
as with the coupled algorithm. 
\remove{Also not shown is a convergence study of 
uniform grids with refinement corresponding to levels 4, 5 and 6, where we have verified
second order convergence at peaks and troughs using a radial reference
solution.}

\begin{figure}[h]
\centerline{\includegraphics[width = .6\textwidth]{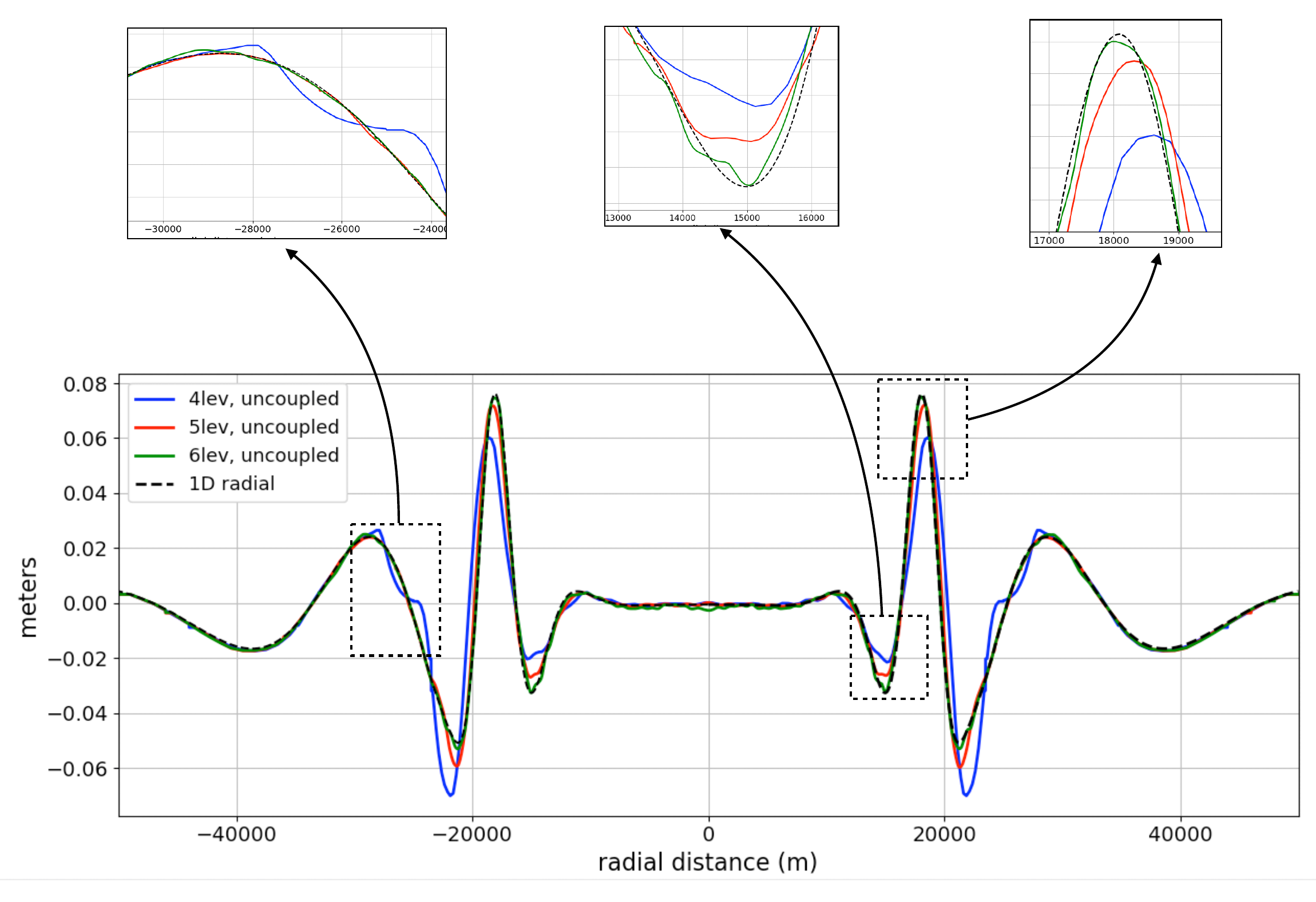}}
\caption{Example showing 4, 5 and 6 level adaptive runs with subcycling
in time converging to the 1D radial reference solution.  
\label{fig:convXsec}}
\end{figure}

In \cref{fig:convXsec} we compare 4, 5 and 6 level adaptive solutions to show
convergence, but 
since the adaptive runs do not refine everywhere and have larger
truncation errors at patch boundaries, it is hard to determine the
convergence rate.
 Note the ``bumps" in the zoom at level 4, which are gone by levels 5
and 6. All runs show small scale oscillations that decrease with refinement level. 
Discrepancies are highlighted in the zooms, but the the results are clearly
converging to the reference 1D radial solution.
There are tiny glitches at patch interfaces 
that can be seen in the zoomed plots.

\FloatBarrier

\section{Computational Examples}\label{sec:ex}

\add{In this section we present computational results for three different
test problems that exercise the new algorithms and software.
We have not tried to optimize the runs in a high-performance computing
environment, but we do provide timings
to demonstrate that it possible to solve
significant problems on a laptop. In all cases the code was run on a
MacBook Pro with the Apple M1 chip using 6
threads for openMP and 6 MPI processes.  The linear systems are solved using
PETSc with a Krylov method and algebraic multigrid as a preconditioner.
We reuse the preconditioner for multiple time steps, which
helps reduce the run time.  
We changed the default convergence criteria
from a relative tolerance of $10^{-5}$ to $10^{-9}$.
Other PETSc parameters can be found in the code repository for this paper
\cite{zenodo}, along with the source code for all three of these examples.}

\subsection{Radially symmetric ocean, shelf, and beach}\label{sec:radocean}
\begin{figure}[t]
\hfil
\includegraphics[width=0.65\textwidth]{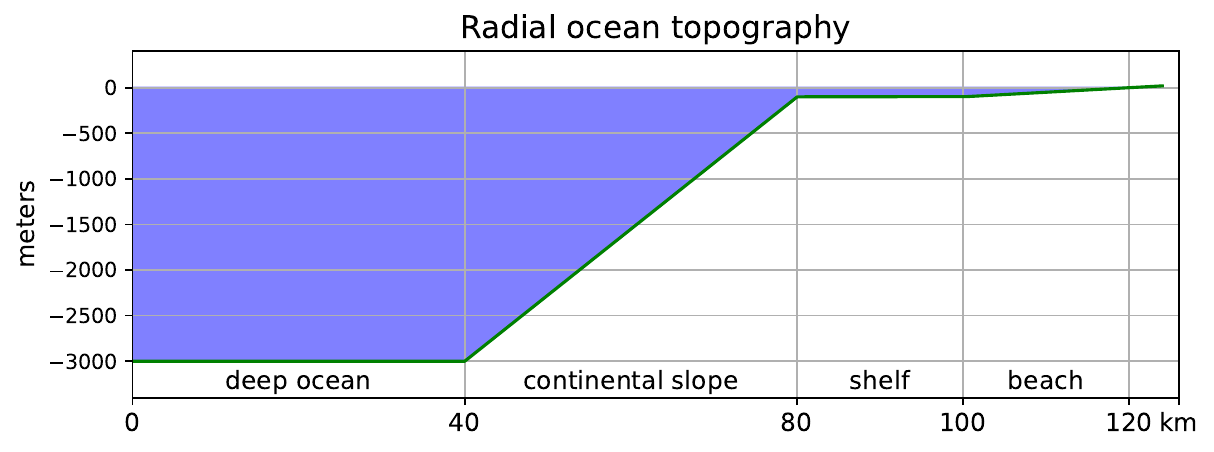}
\hfil
\caption{\label{fig:radocean_topo} Topography for the radially symmetric
ocean test case.
}
\end{figure}

We test both the accuracy and stability of our AMR implementation for a
two dimensional problem where varying topography, dispersion, and nonlinearity
are all important.  We construct a radially symmetric
``ocean'' bounded by a continental shelf followed by a beach. The radial
profile of the topography is shown in \cref{fig:radocean_topo}.
The ocean has a flat bottom with 3000 m depth out to radius of 40 km, a
continental slope from there to 80 km followed by a flat shelf with depth
100 m. Starting at 100 km a beach with slope 1:200 starts to rise and the
initial shoreline is at a radius of 120 km.  As initial data we take a
stationary Gaussian hump of water with surface elevation
\begin{equation}\label{etarad}
\eta(x,y,0) = 20 \exp(-(r/10000)^2)
\end{equation} 
so that the amplitude is 20 m and it decays over roughly 10 km.
We use the one-dimensional radially symmetric SGN equations \eqn{sgnradial}
to compute a reference solution on this topography. The
one-dimensional version of GeoClaw used for this does not support AMR but
uses a finite volume grid with variable spacing chosen so that the Courant
number is close to 1 everywhere that the depth is greater than 50 m,
transitioning to
equal grid spacing in shallower water and onshore.  A total of 10,000 grid
cells were used from $r=0$ to 124 km, giving a grid spacing that varied from
roughly 80 m in the deep ocean to 5 m near shore. The standard GeoClaw
wetting-and-drying algorithm is used on shore. With the initial
conditions used here, the inundation proceeds onshore about 1.2 km, with a runup
elevation of roughly 6 m.  The SGN equations are used in water deeper than 5 m,
switching to SWE in shallower water and onshore for the inundation modeling.

For the 2D simulation we used topography obtained by rotating the 1D profile.
Five levels of refinement
were used with mesh spacing ranging from 2 km on level 1 to 5 m on level 5.
The refinement ratios \add{going from coarsest level 1 to finest level 5 were}
4,5,10,2.  The problem was solved on one quadrant of the
full 2D domain, for $0 \leq x,y \leq 126$ km, with 
solid wall boundary conditions at $x=0$ and $y=0$. In order to test
the stability and accuracy of AMR near grid interfaces, we focused the
refinement on the waves propagating along the diagonal where $x=y$, which 
hit the shore at $x=y=\sqrt{2} \times 120 \approx 169.7$ km. 
Levels 3--5 were allowed around the shelf and beach region.
\remove{ and were further }\add{The cells flagged for refinement were also}
restricted so the finest level patches were created only near the
diagonal, since this is the region of interest for this comparison.

\begin{figure}[th]
\hfil \includegraphics[width=0.95\textwidth]{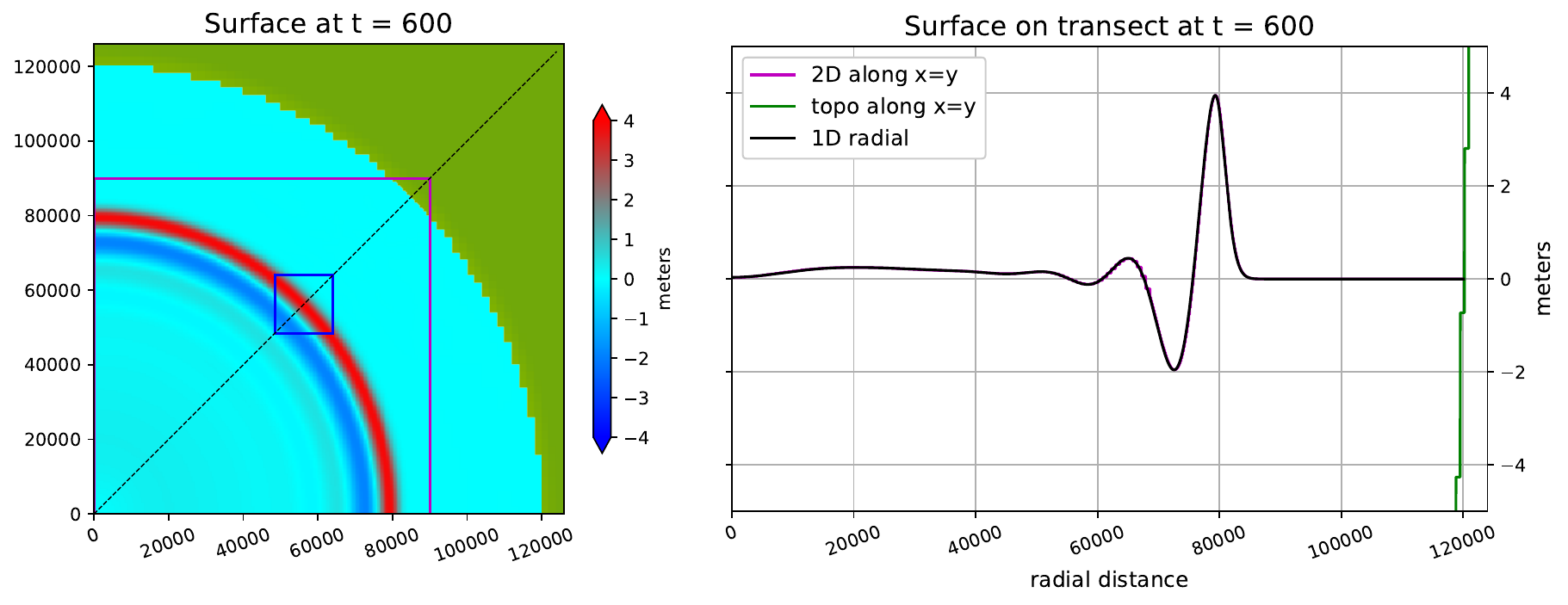} \hfil
\vskip 5pt
\hfil \includegraphics[width=0.95\textwidth]{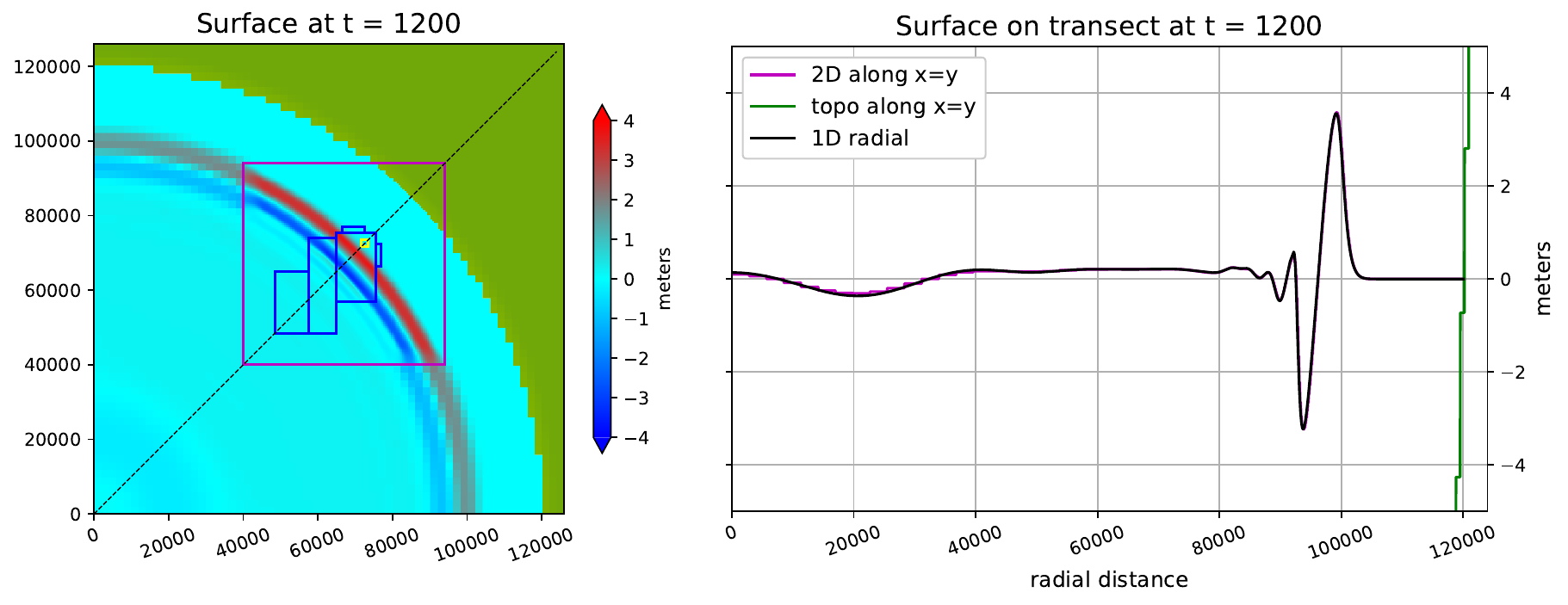} \hfil
\vskip 5pt
\hfil \includegraphics[width=0.95\textwidth]{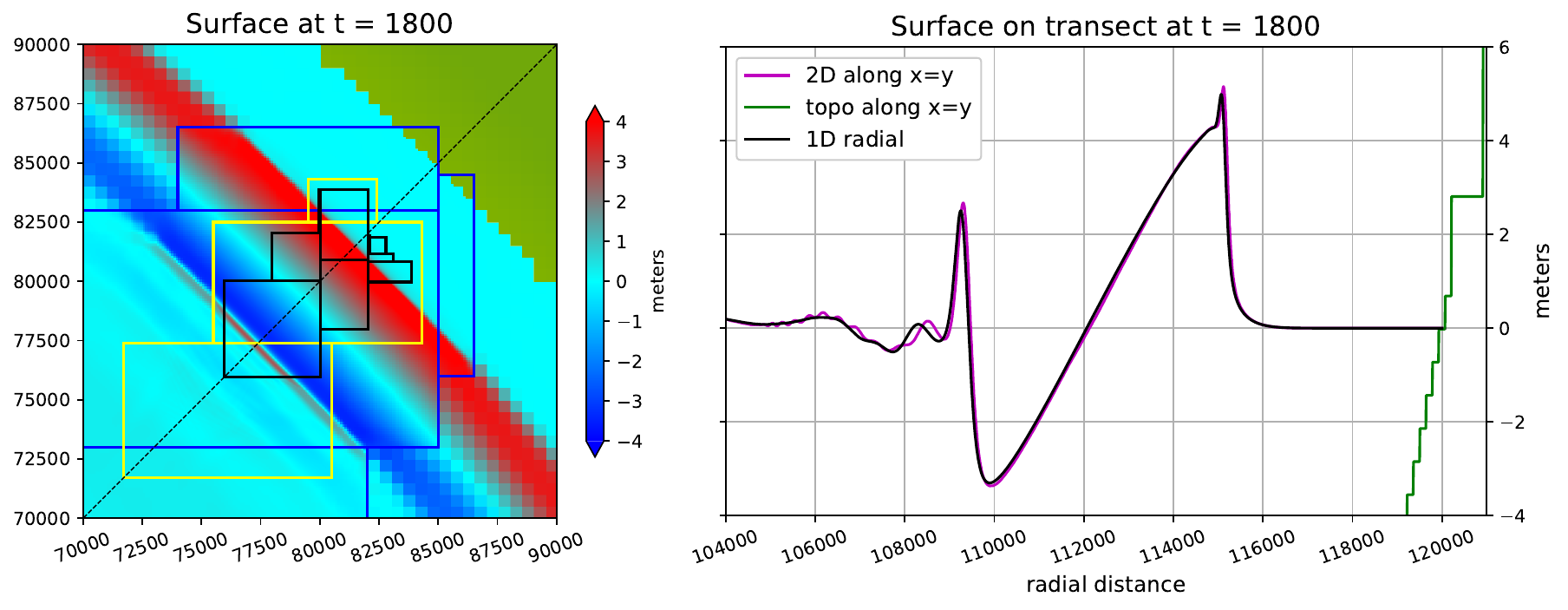} \hfil

\caption{\label{fig:radocean1} Solution to the radially symmetric ocean test
problem at several times.  The left column shows a colormap of surface
elevation $\eta(x,y,t)$ over the positive quadrant.
Colored rectangles show grid patches 
(\add{level 2: magenta,}
level 3: blue, level 4: yellow, level 5: black) and the dashed line is the $x=y$
diagonal. 
The right column shows a transect of the 2D
solution along the diagonal at each time, along with the solution to the 
1D radially symmetric equations at the same time.
Plots are zoomed in near the shore at time $t=1800$ seconds.
}
\end{figure}

\begin{figure}[t]
\hfil \includegraphics[width=0.95\textwidth]{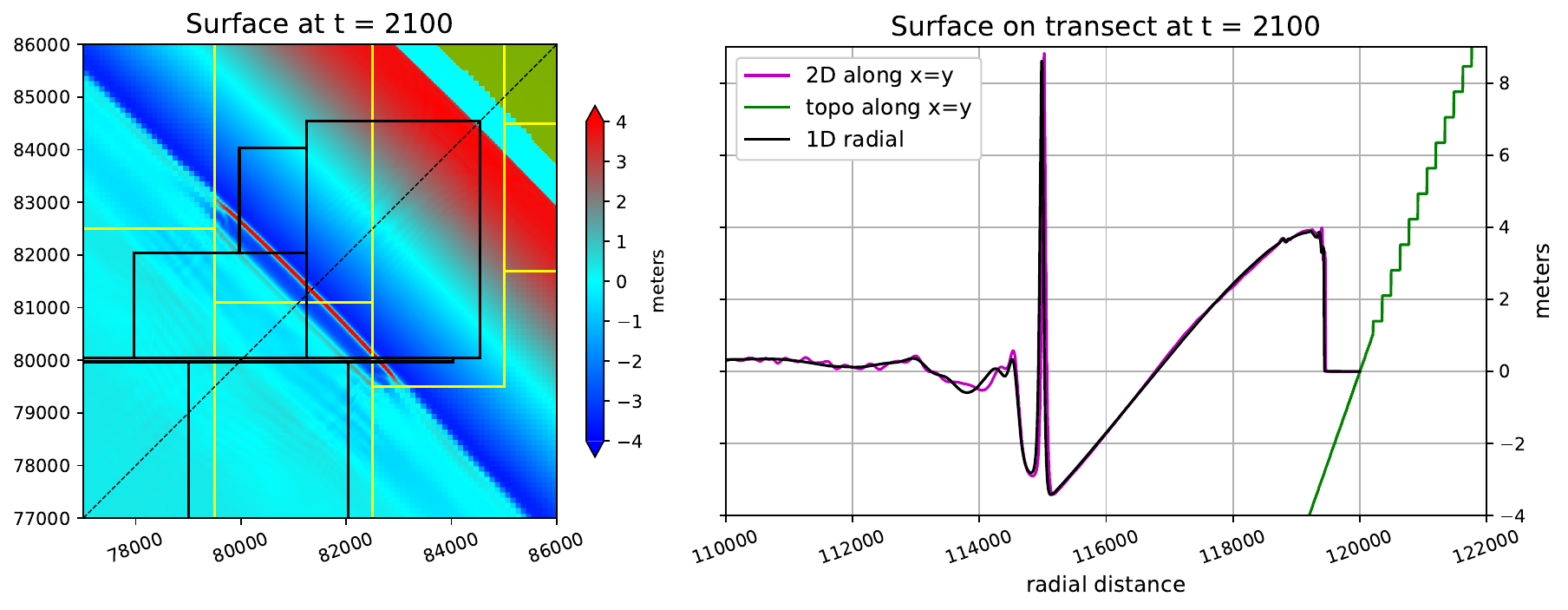} \hfil
\vskip 5pt
\hfil \includegraphics[width=0.95\textwidth]{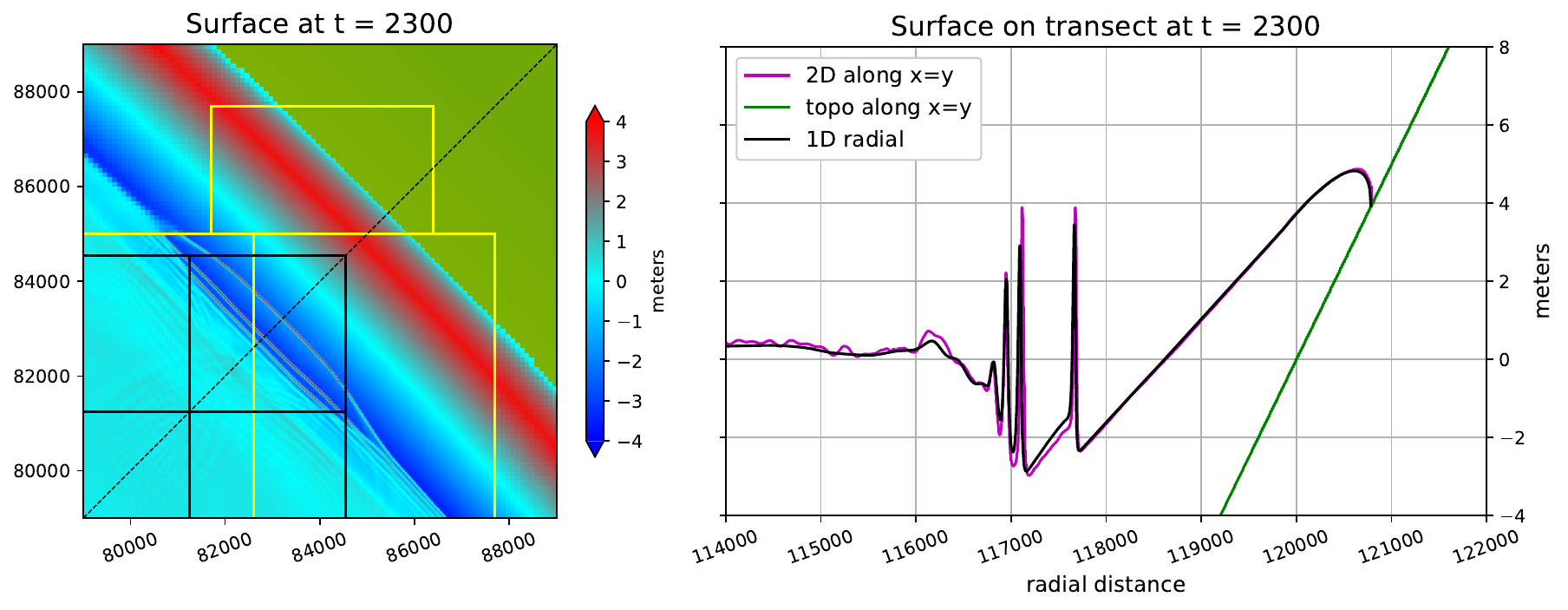} \hfil
\vskip 5pt
\hfil \includegraphics[width=0.95\textwidth]{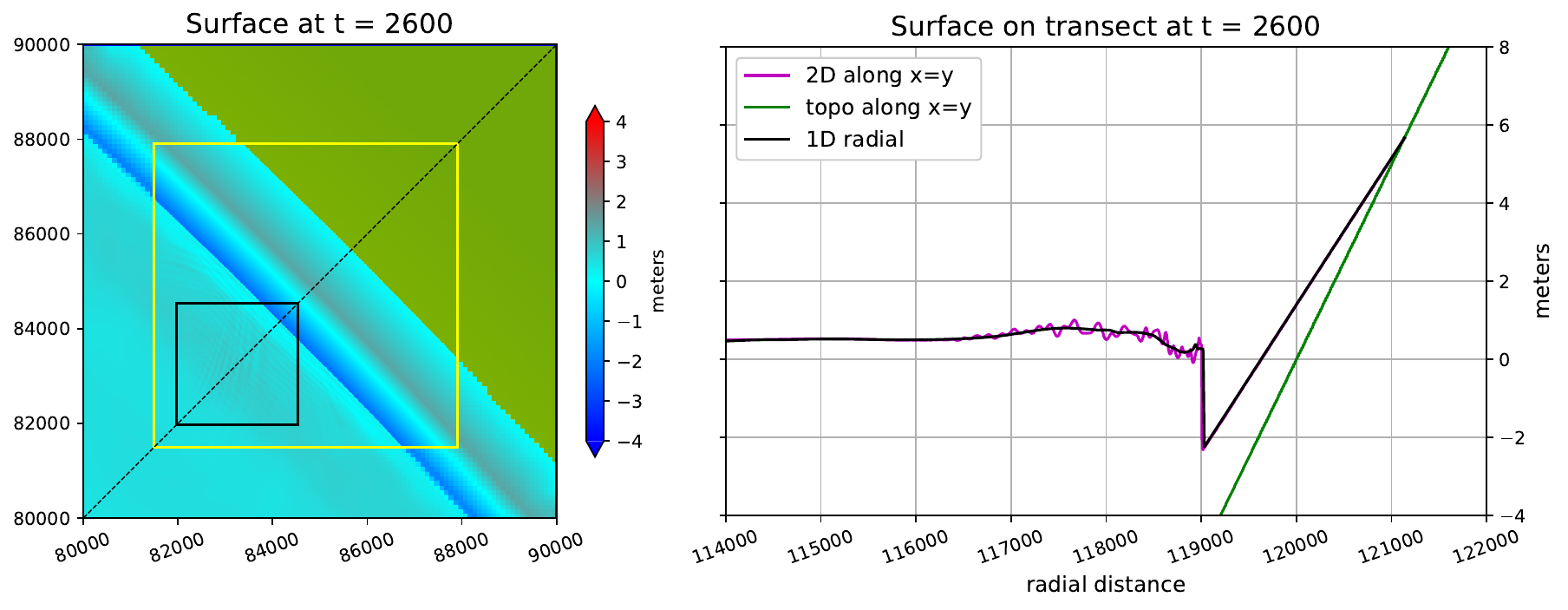} \hfil

\caption{\label{fig:radocean2} Continuation of \cref{fig:radocean1}
at later times (zoomed in near the shore).
\add{Note that the axes change in each plot to better show the signficant
features at each time.}
}
\end{figure}

\Cref{fig:radocean1,fig:radocean2} show the computed
solutions at several times, both the 2D surface elevation over the
quadrant and the transect along the diagonal. The latter is
plotted together with the computed solution to the 1D radial equations to
assess the accuracy.  
Note that by $t=600$ seconds the dispersion has created an oscillatory wave
train that is well captured on level 3. By time $t=1200$, shoaling on the
continental shelf causes the wavelength to decrease, and then by $t=1800$
the nonlinearity in the shallower water has caused the
waves to steepen. Rather than breaking,
the dispersion leads to ``soliton fission'', the
appearance of solitary waves that break off from the steepening main waves.
This phenomenon is frequently observed as real tsunamis approach shore.
At $t=1800$ a soliton is just forming near the main peak. These waves are
well captured. At $t=2100$, the leading wave is in water with
depth less than 5 m, where the transition from SGN to SWE occurs, and a
shock (hydraulic jump) has formed.  In deeper water, at roughly $r=115$ km, a
sharply peaked soliton has formed.  Note that this wave is still 100 m wide
and less than 10 m high, so it is reasonable that this wave has not yet broken.
Also note that the 2D solution agrees very well
with the 1D solution in this region
where the grid is refined to a similar resolution, in spite of the fact that
the 2D wave is moving diagonal to the Cartesian grid here. 
At $t=2300$ the wave is moving onshore and has inundated about
800 m inland. The agreement with the 1D radial solution is excellent
onshore and quite good in the region of soliton fission.
Finally, at $t=2600$ the wave has roughly reached the level of maximum runup
and is starting to retreat.

In the 2D plots the accuracy decreases dramatically outside of
the finer grid patches.  Inside the fine grid patches, very
small oscillations can be seen
radiating from the patch boundaries that also appear in the transect plots,
but these oscillations have little effect on the primary waves in the region
of interest.

\clearpage  

The simulation up to $t=2300$ required about \remove{18}\add{6}
hours (wall time) on a laptop as described at the beginning of \cref{sec:ex}.  
\add{(For better visibility of both the underlying solution and the grid patches in \cref{fig:radocean1,fig:radocean2},
we allowed large grid patches (up to $1000\times 1000$ instead of
$100 \times 100$), reducing the
effectiveness of using OpenMP for the explicit shallow water time stepping
and slightly increasing the run time).}  
A simulation using the shallow 
water equations with the same refinement levels and strategy required about
30 minutes of wall time.  However, the SWE solution (not shown)
does not contain the narrow solitons seen
in the SGN solution, and so it would not require the same level of
refinement to obtain an accurate SWE solution for this problem, potentially
further reducing the time required for SWE.


\subsection{Japan 2011 Earthquake Tsunami}\label{sec:tohoku}
The Great East Japan (Tohoku) Earthquake of March 11, 2011 was studied by
Popinet \cite{Popinet2012} using the shallow water equations.
He further used it as a test problem for the dispersive SGN equations
in \cite{Popinet:2015}, providing some comparisons between the two
sets of equations.  We present a similar comparison for our AMR
implementation of the same version of the SGN equations so that we can compare
with his figures. Because
the primary tsunami wavelengths generated by this megathrust event
were so long, the addition of dispersive terms has relatively little
effect on the solution for the most part.  This is generally true
when modeling tsunamis generated by large subduction zone earthquakes,
the source of most devastating historical tsunamis, and the use of
the shallow water equations without dispersive terms is often well
justified \cite{Glimsdal2013,kirby_dispersive_2013}.
However      
there are some higher frequency waves introduced by the dispersion that
match some of the observations slightly better.
Below we also use this example to investigate a conjecture
made in \cite{Arcos2015}
concerning the arrival times of the peak waves. 

This problem is also a good
test of the stability of our AMR implementation of the SGN equations on
real topography from transoceanic to harbor scale, with realistic refinement
patterns and running conditions similar to those often used in practice.
Six levels of adaptive refinement were used for the simulations presented
below, starting with a level 1 grid having 2 degree resolution (220 km) and
using refinement ratios 5,6,4,6,30 at successive levels to reach 1/3
arcsecond (10 m) resolution in Kahului Harbor on the island of Maui, HI.
In the deep ocean at most level 4 was used (1 arcminute), with frequent
regridding to follow the leading wave traveling to Hawaii without
over-resolving other regions.  
We used identical running conditions as in \cite{Arcos2015} with one
exception: there the deep ocean was refined only to level 3 (4 arcminutes)
but better resolution is required to model the more oscillatory dispersive
waves, and even to capture the SWE waves properly close to the source.
In these new computations, 4 levels were allowed in the deep ocean for
both the SGN and SWE simulations.

Several different reconstructions of the tsunami source
(seafloor motion) from the 2011 earthquake have 
been developed by different groups. 
Popinet used a model developed by the UCSB group of
Shao et al.\ \cite{shao_focal_2011}.
For comparison we present some results with this UCSB model.
Then we also present some results using a source model developed
by Fujii et al.\ \cite{fujii_tsunami_2011},
which was found in \cite{MacInnes2012} to be one of the best at
replicating nearfield run-up and DART buoy time series in a
comparison of GeoClaw results for ten proposed sources for this event
(which also included the UCSB source).

\begin{figure}[t]
\hfil
\includegraphics[width=0.45\textwidth]{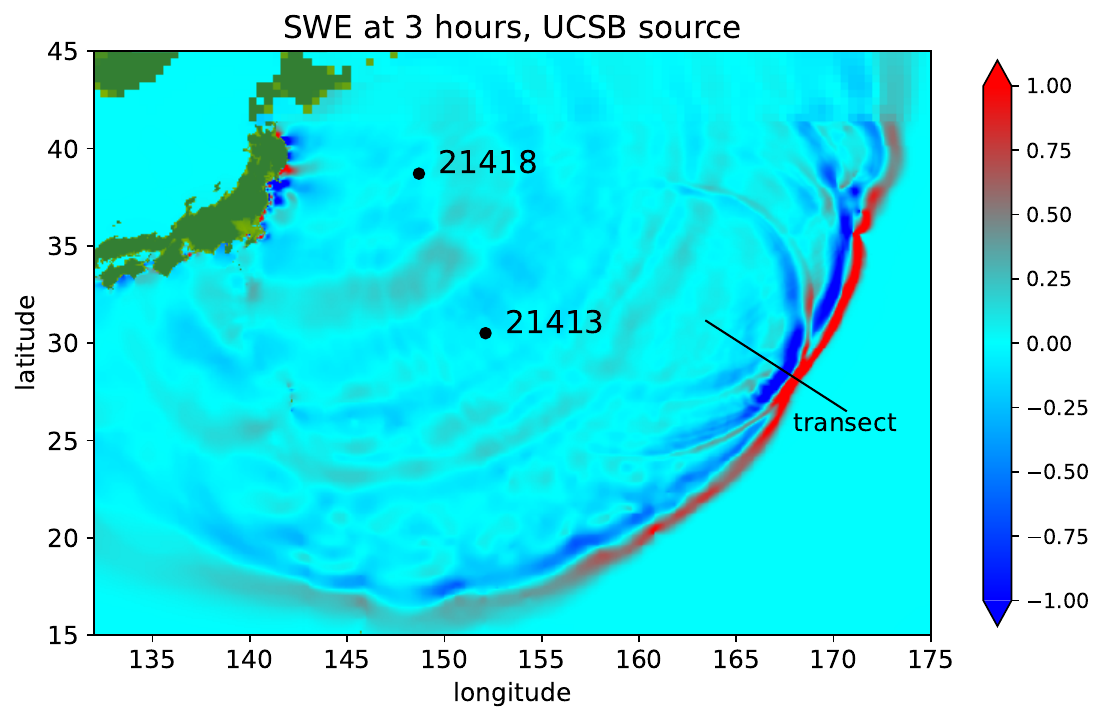}
\hskip 5pt
\includegraphics[width=0.45\textwidth]{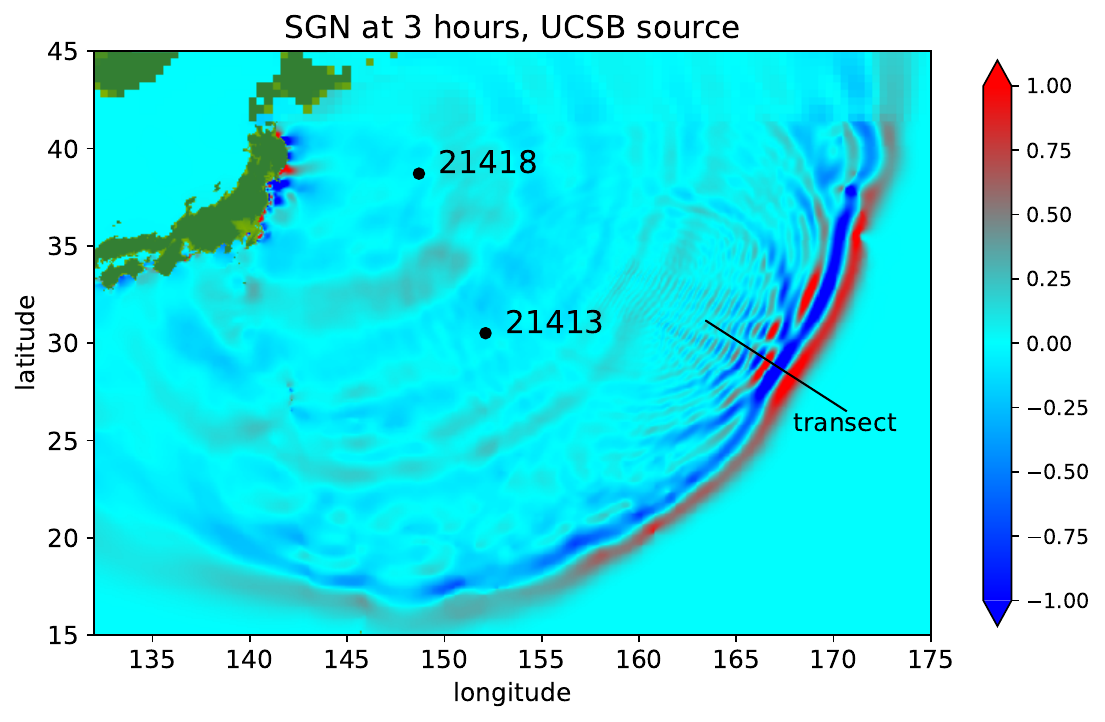}
\hfil
\caption{\label{fig:tohoku3hrs} Simulations of the 2011 Tohoku tsunami 3
hours after the earthquake, using the UCSB source model, using SWE (left)
or SGN (right). Colors saturate at $\pm 0.5$ m relative to sea level in order
to accentuate the small oscillatory waves following the main peak in the
dispersive SGN simulation. The DART 
gauge locations are also indicated, for which time series data is
plotted in \cref{fig:tohokuDART}. Compare this plot to the bottom row
of Figures~12 and 13 in Popinet \cite{Popinet:2015}.
The black lines show the transect along which the solution is plotted
in \cref{fig:tohokuTransect3hrs}.
}

\vskip 10pt 

\hfil
\includegraphics[width=0.70\textwidth]{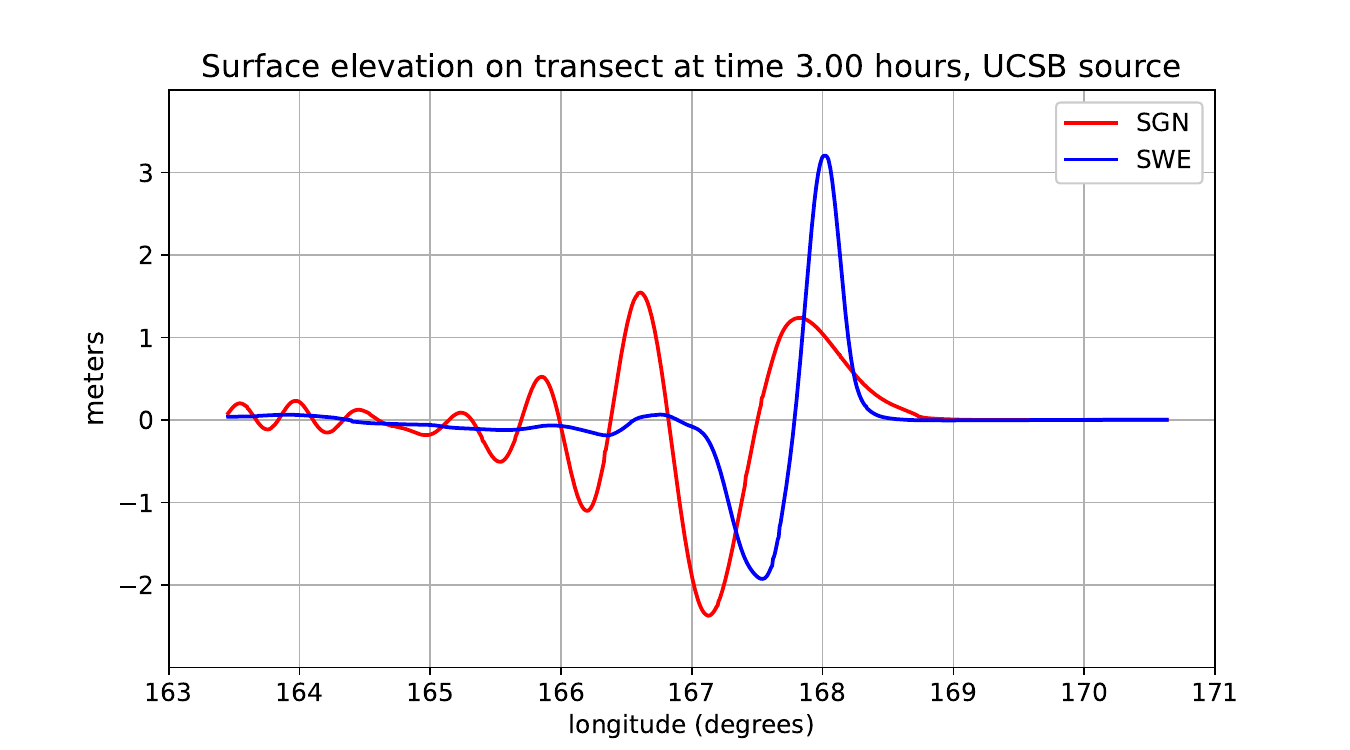}
\hfil
\caption{\label{fig:tohokuTransect3hrs} Cross section of the solution from
\cref{fig:tohoku3hrs} at 3 hours, along the transect shown in that figure,
from (163.454E,~31.1458N) to (170.631E, 26.5282N) as provided by S.~Popinet
(personal communication).
Compare this plot to Figure~15 in Popinet \cite{Popinet:2015}.
}
\end{figure}

\Cref{fig:tohoku3hrs} shows the tsunami generated by the UCSB source at
3 hours, as computed both with SWE and SGN. Differences are relatively
minor, but the main outgoing wave is seen to have some dispersive ripples in
the SGN calculation that are not present with SWE.
The differences are more clearly seen when a transect of the surface elevation
at a fixed latitude is plotted, see \cref{fig:tohokuTransect3hrs}.
This is similar to Figure~15 in \cite{Popinet:2015}.
It should also be noted that the Basilisk \cite{basilisk} model used by
Popinet does cell-by-cell refinement on a quadtree structure, and that
he refined in
somewhat different regions than in our patch-based GeoClaw simulations,
so direct comparison is not entirely possible.
 
\begin{figure}[t]
\hfil
\includegraphics[width=0.75\textwidth]{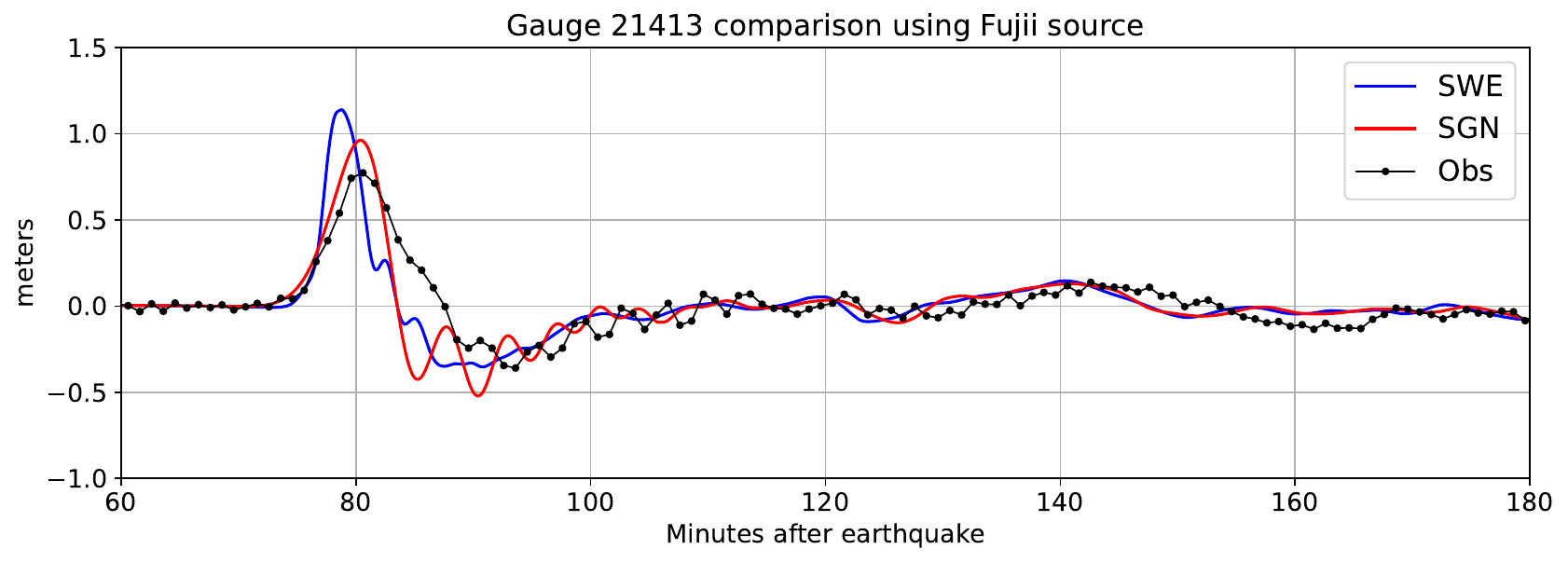}
\hfil
\vskip 5pt
\hfil
\includegraphics[width=0.75\textwidth]{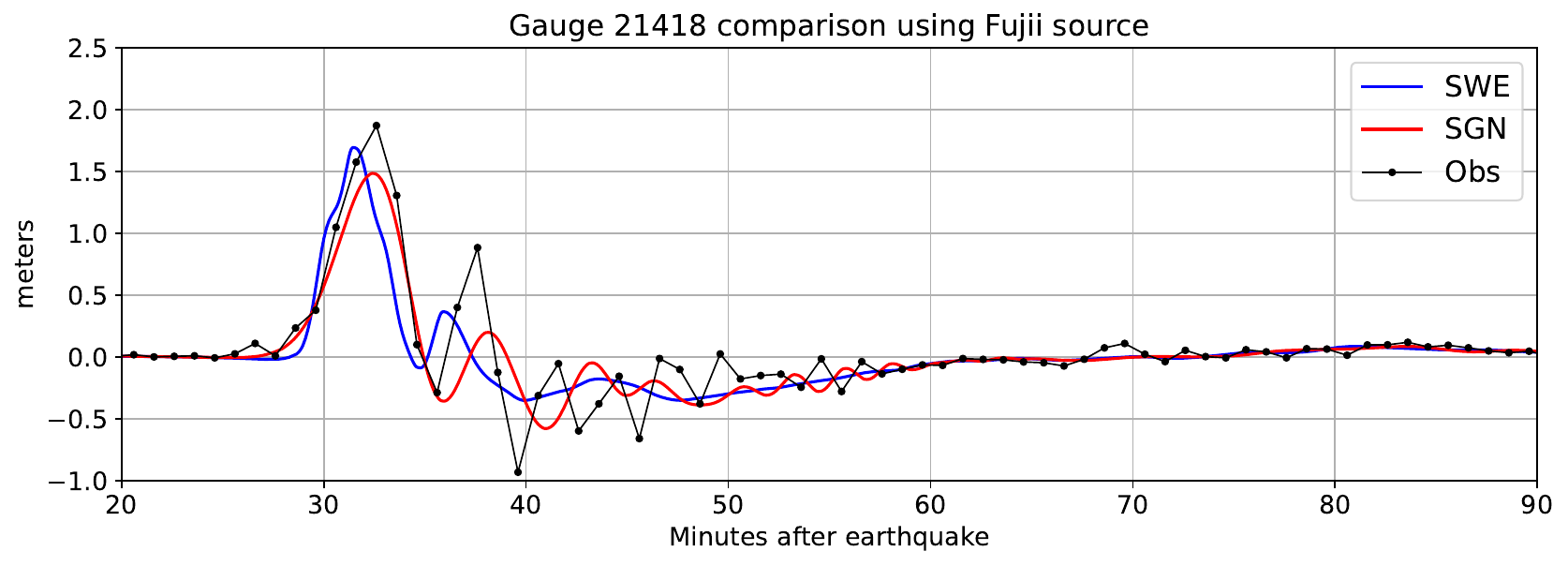}
\caption{\label{fig:tohokuDART} DART gauge time series at two locations
marked in \cref{fig:tohoku3hrs}.
Detided observations (Obs) from March 11, 2001 
are compared with SWE and SGN computations using the Fujii tsunami source.
The SGN results show additional oscillations that appear to match the
period but not phase of the under-resolved observations.
}
\end{figure}

\Cref{fig:tohoku3hrs} also shows the location of two DART sensors that
recorded the tsunami passing by.  \Cref{fig:tohokuDART} shows the
de-tided DART data observed at these locations (from the data archived with 
\cite{MacInnes2012}), together with the time series computed with both the
SWE and SGN versions of GeoClaw. Here we show results using the Fujii source,
which in general matches the DART observations better.
As expected, dispersion makes relatively little
difference for these long wavelength waves, but does add some oscillations
at roughly the same period as seen in the observed data, 
although the latter is at a
sampling rate of 60 seconds and is not well resolved at shorter wavelengths.

The Fujii tsunami source of \cite{fujii_tsunami_2011}
was used by Arcos and LeVeque \cite{Arcos2015}
to model tide gauge observations and tsunami current speeds at a number
of points in Hawaii.  In that paper the SWE were used, and many of the computed
time series agreed remarkably well with detided observations after shifting the
results in time by about 10 minutes.  Several possible reasons for this time
shift were discussed in \cite{Arcos2015}, including the fact that an
instantaneous seafloor displacement was used in the model.  The possible effect
of dispersion on the time shift was roughly estimated
using the average ocean depth and wavelength of the leading wave.  
It was estimated that the arrival time might be shifted by dispersion 
by about 1\%, less than 5 minutes in Hawaii, where the first wave
arrives roughly 460 minutes after the earthquake.

\begin{figure}[h!]
\hfil
\includegraphics[width=0.75\textwidth]{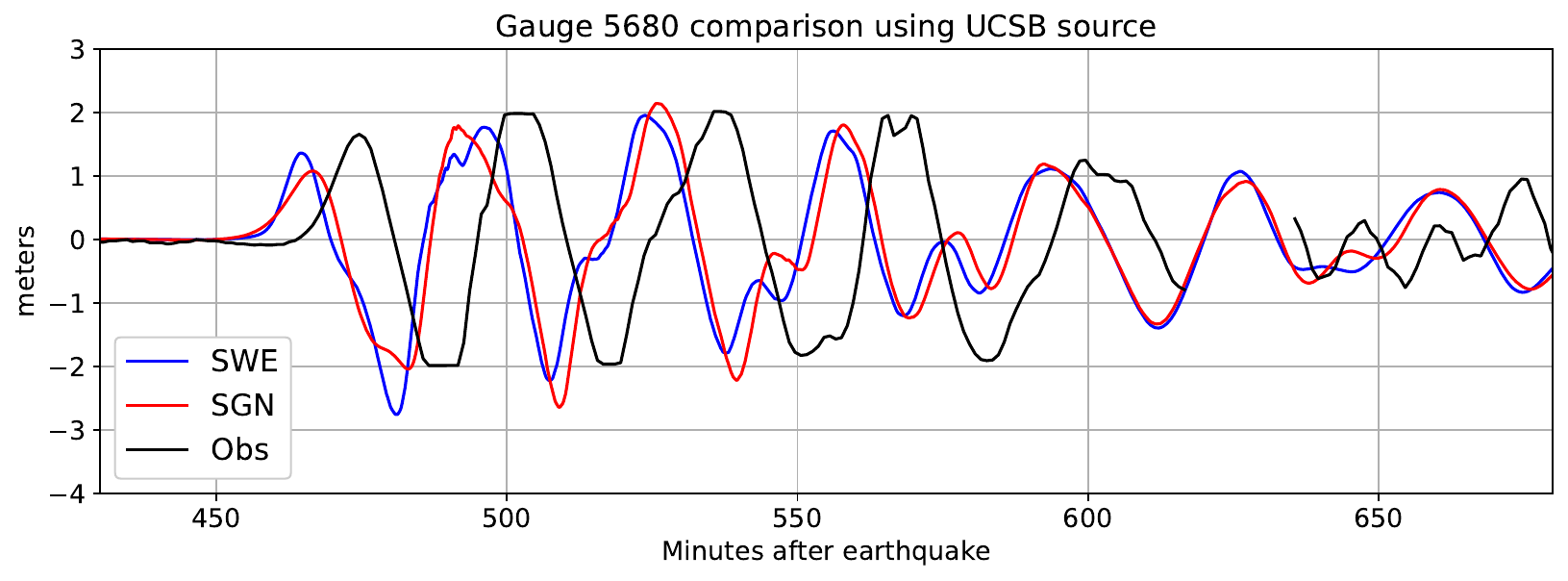}
\hfil
\vskip 3pt
\hfil
\includegraphics[width=0.75\textwidth]{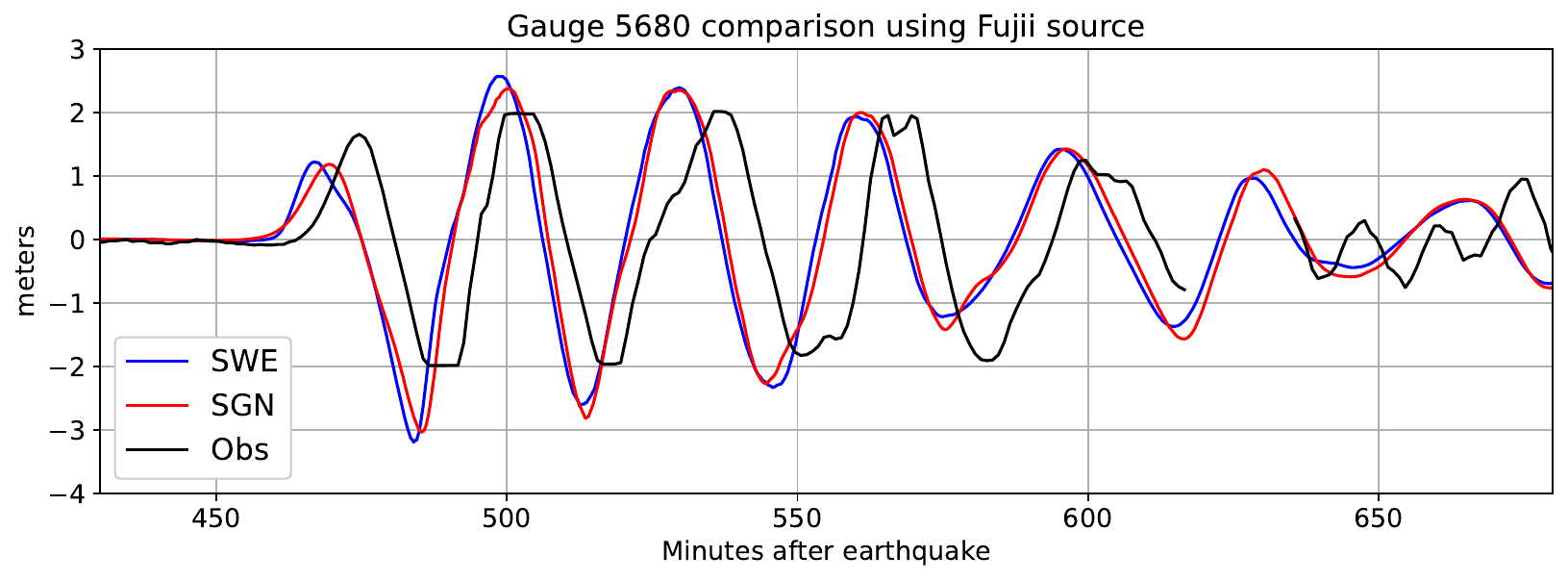}
\caption{\label{fig:tohokuTG} Tide gauge time series in Kahului Harbor.
Detided observations (Obs) are compared with SWE and SGN computations
using the UCSB (top) and Fujii (bottom) tsunami sources.
Note that the physical tide gauge bottoms out at $-2$ m while the synthetic
gauge is in slightly deeper water, and some observational data is missing
after 600 minutes.}
\end{figure}

\Cref{fig:tohokuTG} shows the detided observation at one tide gauge
located in Kahului Harbor (Maui), along with the computed time series
using both the SWE and SGN equations, with no time shift. 
The top figure shows the computational results obtained with the UCSB source,
while the bottom figure shows those obtained
with the Fujii source, for comparison to \cite{Arcos2015}.  In the Fujii
plot, the computed waveforms look nearly identical except for a small time
shift. With SWE,
the first maximum tsunami amplitude occurs about 8 minutes early compared to
the observation.
Using SGN, the peak appears roughly 3 minutes later than
with SWE, so only 5 minutes early relative to the observations.
This is consistent with other findings in the
literature that even with the use of dispersive equations, numerical
simulations show waves arriving about 5 minutes earlier than the
observations in Hawaii for this event (e.g.
\cite{kirby_dispersive_2013,cheung_surges_2013}).
The time series computed with the UCSB source show similar shifting of the
first peak, by about 2 minutes. We also note that the differences
observed between the SWE and SGN results in each case are
small relative to the differences between using the two different source
models, neither of which is an exact representation of the actual earthquake.
\remove{This again confirms that the use of dispersive
equations is relatively unimportant for such events.}

These calculations were performed on a laptop as described
at the beginning of \cref{sec:ex}.
Using the UCSB source, the
shallow water run to 12 hours of simulated time took
roughly \remove{100}\add{23} minutes of wall time, while the SGN run required
\remove{366}\add{145} minutes.
\remove{Using the refinement criterion we chose, a smaller portion of the ocean was
refined when using the Fujii source. As a result, run times were less for both
SGN and SWE, but they exhibited a similar factor of 3
difference between the SWE and SGN runs.}


\subsection{Hypothetical Asteroid Impact}\label{sec:expdc}
Next we simulate the tsunami that results from a hypothetical
asteroid impact and subsequent crater in the Pacific ocean
off the coast of Washington. The paper \cite{ICM22} contains more details and
sample simulation results for this example using the MS equations.
Here we use the SGN equations and we will compare to results in \cite{ICM22}, 
since there no ``data'' for this experiment.
A similar test problem was
also used in \cite{NASA-AGT-2018} for a simulation using the non-dispersive
SWE for a much smaller asteroid.
That paper also highlighted the need for dispersive simulations
to properly model propagation and possibly runup of asteroid-generated tsunamis.

\begin{figure}[h]
\hspace*{-.4in}
\mbox{
\includegraphics[height=1.7in]{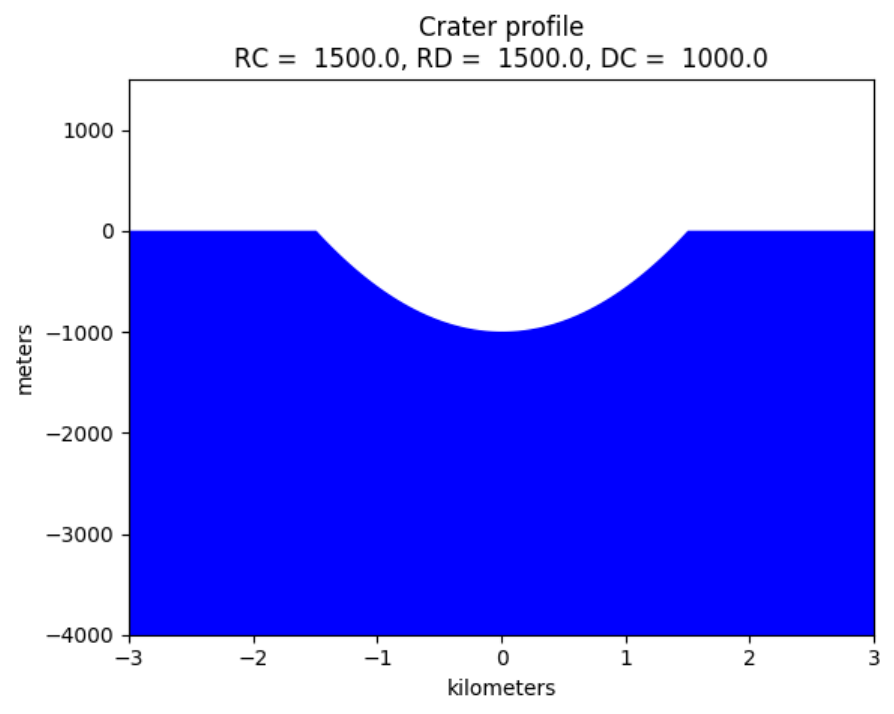}
\includegraphics[height=1.7in]{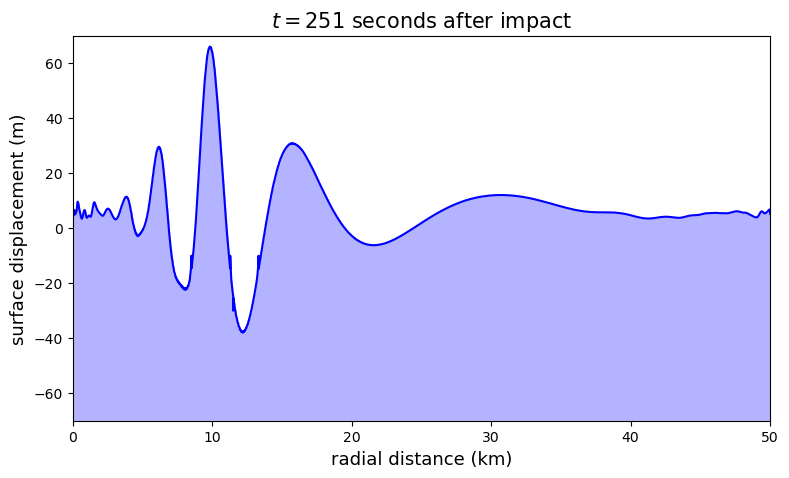}
\includegraphics[height=1.7in]{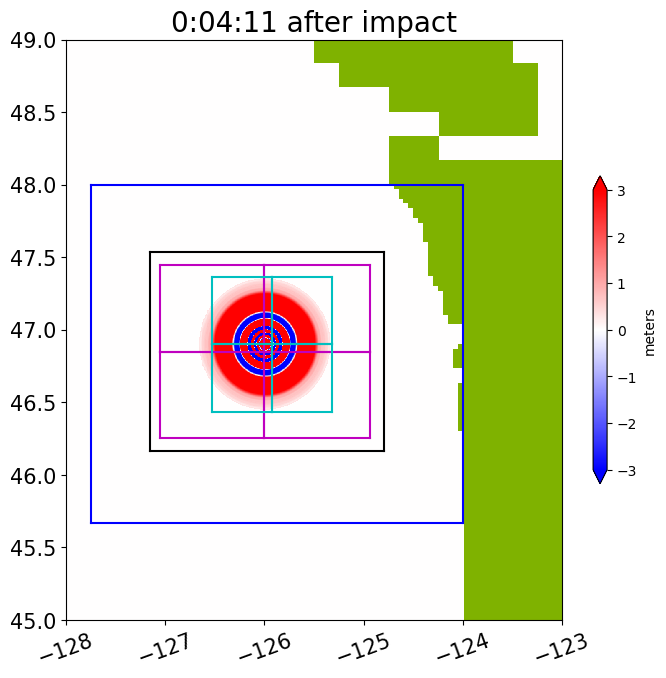}
}
\caption{(a) Initial conditions of a static crater of depth 1
km.~and diameter 3 km., which were the initial conditions for the hydrocode
simulation. (b) The radially symmetric results of the hydrocode
simulation at 251 seconds, used to initialize the GeoClaw Boussinesq
simulation. (c) Radially symmetric GeoClaw initial conditions and mesh,
with the initial grid patches outlined.  \label{fig:craterfig}}
\end{figure}

We use the same problem setup as in \cite{ICM22}. 
The crater is 1 kilometer deep with a diameter of 3 kilometers, as
depicted in \cref{fig:craterfig}a.
The crater was placed 
approximately 150 km west of Grays Harbor, a well-studied area
due to the potential for Cascadia subduction zone earthquakes. 
The initial tsunami surface elevation, shown in \cref{fig:craterfig}b,
is taken after 251 seconds of a radially symmetric hydrocode
simulation starting with a static crater. The velocity is initialized
using an eigenvector from the linearized shallow water equations to
give an approximately radially outgoing wave.
The simulation uses 7 levels of mesh refinement with
refinement factors \add{from coarsest to finest}
of 5,3,2,2,5,6. Initially there are 5 levels, as shown
in \cref{fig:craterfig}c. The last two levels are added as the wave
shoals on the continental shelf and then approaches shore.  The
coarsest level has $\Delta y$ = 10 arcminutes, and the level 7 grids have
$\Delta y$ = 1/3 arcsecond. On each grid $\Delta x = 1.5\Delta y$
so that the finest-level computational grids are at a 
resolution of roughly 10 m in both $x$ and $y$.
The equations switch from SGN to SWE where the initial water depth
is less than 10 meters, regardless of wave height.
(We observed no difference in inundation when using 5 or 2 meters
for the threshold).

\begin{figure}[t]
\centering
\hfil
\includegraphics[width=.45\textwidth]{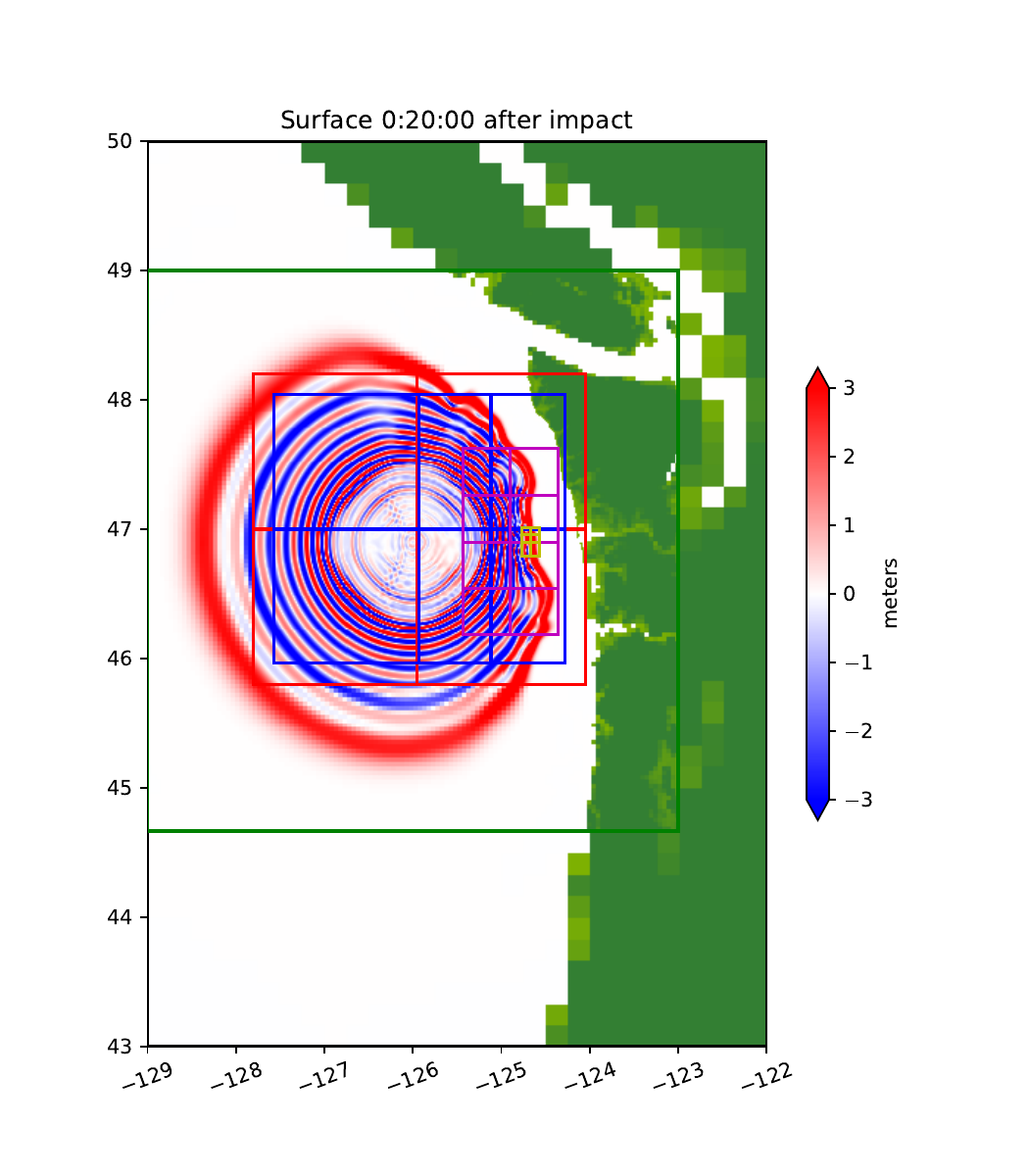}
\hspace*{.25in}
\includegraphics[width=.45\textwidth]{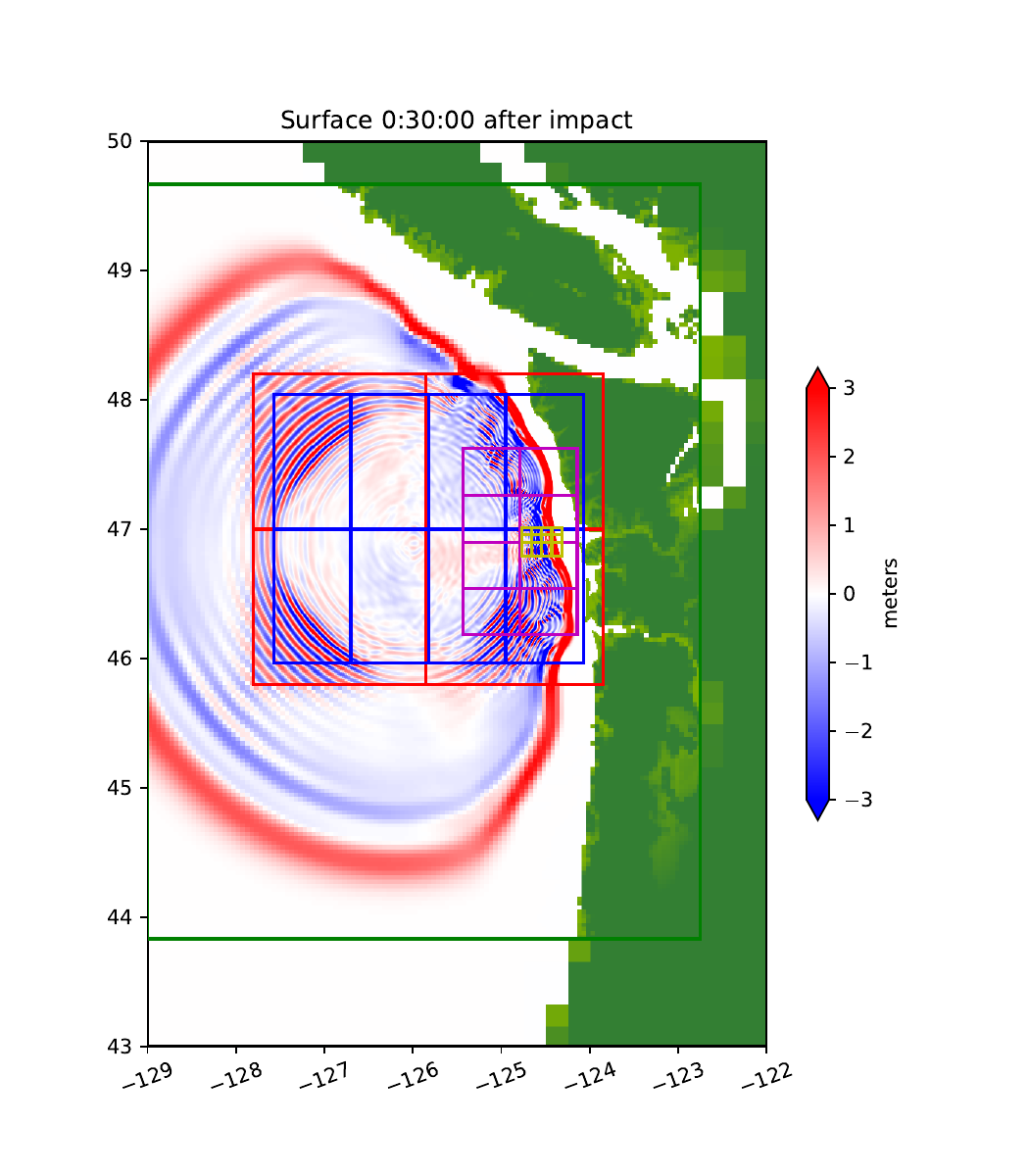}
\hfil
\caption{Surface elevation computed using SGN at the indicated
times. The outlines of the rectangles show the grid patch boundaries, and
are color-coded by level of refinement as described in the text.
\label{fig:waveHeights}}
\end{figure}

\begin{figure}[h!]
\centering
\hspace*{-.10in}
\mbox{
\includegraphics[width=.49\textwidth,trim=0 0 36 0,clip]{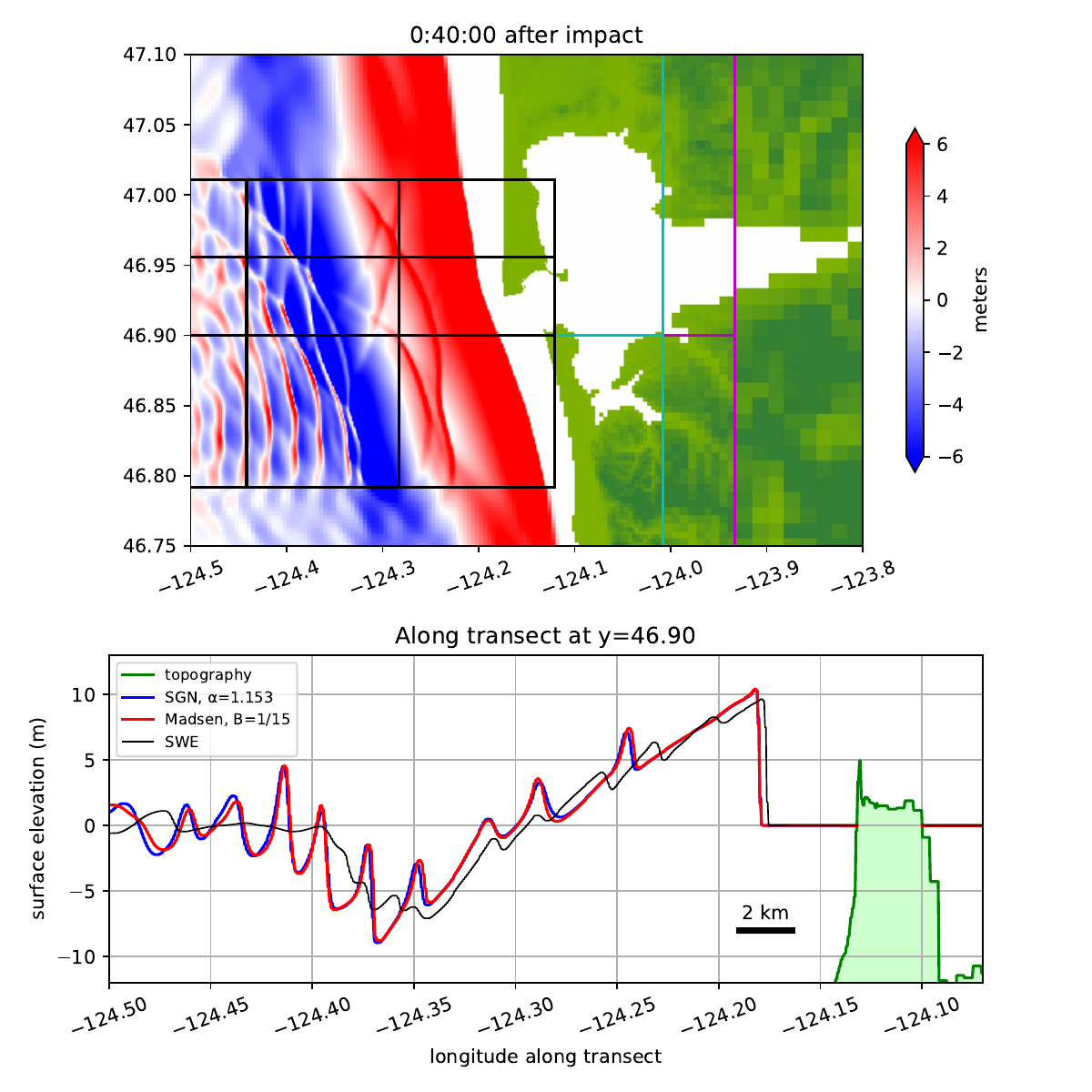}
\includegraphics[width=.49\textwidth,trim=25 0 25 0,clip]{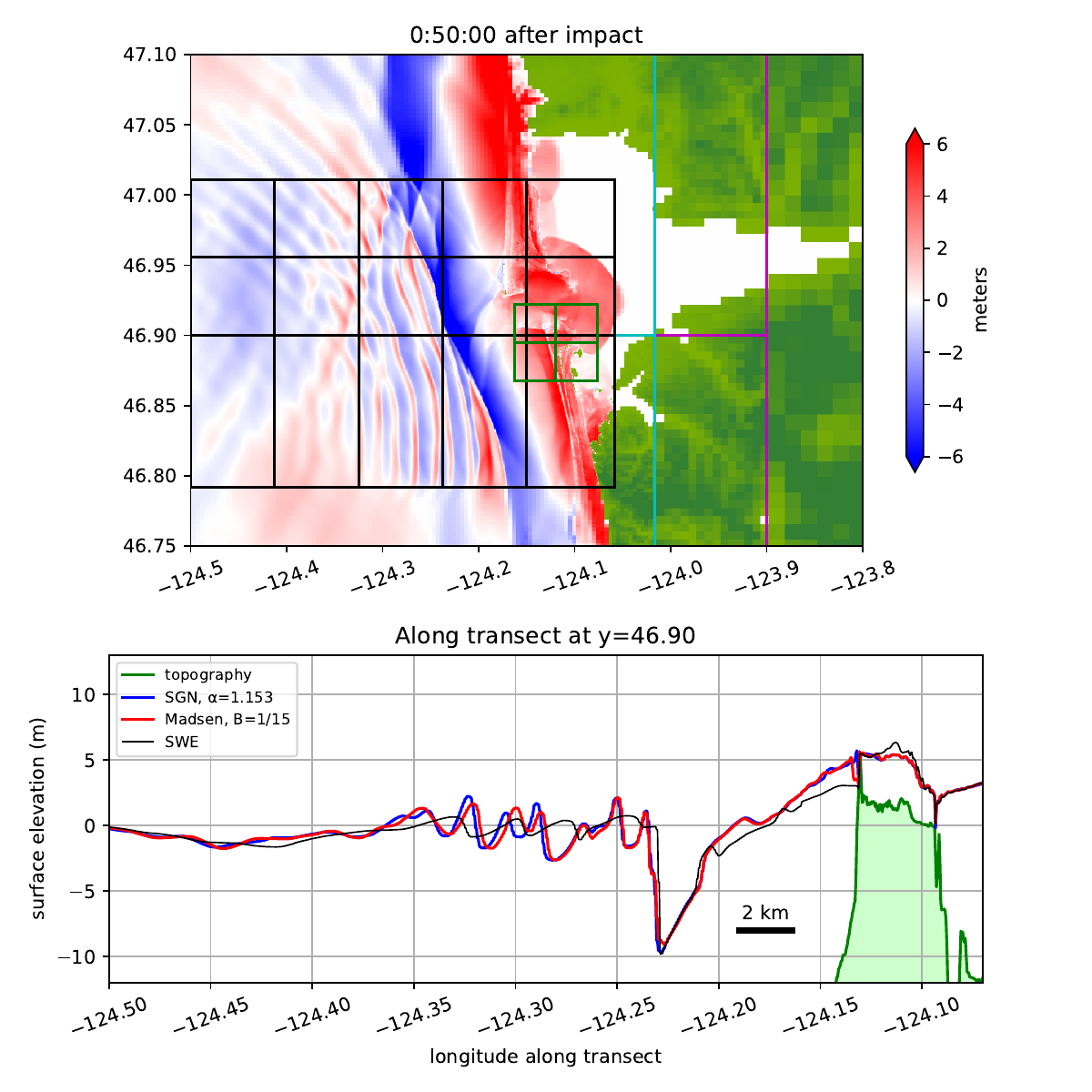}
}
\caption{Top: Surface elevation computed using SGN at the indicated
times. The outlines of the rectangles show the grid patch boundaries, and
are color-coded by level of refinement as described in the text. Bottom: Surface elevation along a
transect at latitude 46.90 cutting through Westport, WA,
highlighting the very small differences between
the MS and SGN solutions. \add{The solution computed with SWE is also shown
for comparison, and does not exhibit the solitary wave formation seen in the
dispersive simulations.} Note that at 50 minutes the Westport peninsula is
refined to one higher level, so the topography is finer in the transect.
\label{fig:solitons}}
\end{figure}

The wave height after 20 and 30 minutes is shown in \cref{fig:waveHeights}.
Both compare extremely well with the comparable figures in \cite{ICM22}. 
There are more grid patches in this
calculation than that paper, since a smaller maximum patch size was
specified to allow for higher parallelism. 
Even the region with soliton fission, shown in \cref{fig:solitons} 
as the waves approach shore, are amazingly similar, so we also show
transects of the solution at a fixed latitude to see
the differences. Since the nonlinear and
dispersive terms in SGN and MS differ, one might have expected
more difference in this area. The refinement to finer levels at later times is
focussed on Grays Harbor, and the plots show that the
two peninsulas at the entrance to the harbor
(containing the communities of Ocean Shores and Westport, WA)
are completely inundated. 
In both figures the patch boundaries are color-coded by level: level 2 is in
blue, 3 is black, 4 is magenta, 5 is cyan, 6 is black and 7 is green.

The simulation up to 50 minutes required about \remove{67}\add{30} minutes 
(wall time) on a laptop as described at the beginning of \cref{sec:ex}. 
A simulation using the shallow
water equations with the same refinement levels and strategy (not shown)
required about \remove{18}\add{4} minutes.

\FloatBarrier

\section{Conclusions and Future work}\label{sec:conc}

We have developed a high-fidelity tool for solving depth-averaged SGN
equations, a form of Boussinesq equations that can be used to
accurately model short-wavelength tsunamis
or other related dispersive wave problems.  Building on the open source
GeoClaw software, patch-based mesh refinement in space and time is implemented
for solving realistic ocean-scale problems.
The major new algorithmic advance is the introduction of an implicit solver
in the AMR code to solve an elliptic equation 
on all grid patches at a given level of refinement.  
The elliptic equation is solved before taking each shallow water time step.
It must also be solved provisionally at the end of each time step on all but
the finest AMR level in order to provide values for space-time interpolation
to ghost cells 
needed to solve the elliptic equation on the finer grid levels,
in exactly the same
way that ghost cell values for the depth and momentum are interpolated.  This
approach allows subcycling in time, so that the equations are advanced on
each level using a time step based on the CFL restriction of the shallow
water equations, allowing optimal efficiency and minimal numerical dissipation.
By using PETSc for solving the elliptic equations, in conjunction with OpenMP
for the explicit hyperbolic steps, an efficient code has been developed in the
GeoClaw framework.
In the near future this will 
be merged into the GeoClaw code base and better documented for 
general use by the scientific community. 
We also plan to apply this software to several modeling problems, in particular
for the NASA Asteroid Threat Assessment Project that originally motivated
and partially funded this work.

In future research, we plan to better explore the transition between SGN
and SWE in the nearshore region in order to better model breaking waves,
which can have an impact on the resulting onshore inundation.
We would also like to reduce the generation of reflections at
interfaces between different refinement levels.  Although these reflections
are generally minor, as seen in the examples we have presented, they do tend
to be larger than those seen when solving SWE.  Better interpolation
procedures may be required for the SGN equations since they involve higher
order spatial derivatives. In addition it may be beneficial to
increase the accuracy of the finite difference discretization of the
dispersive terms. 

We also plan to compare the Boussinesq approach to modeling dispersion,
which requires implicit algorithms, with hyperbolic formulations
(e.g. \cite{bassi_hyperbolic_2020,guermond_hyperbolic_2022}) 
that model dispersion by introducing relaxation source terms.
This requires somewhat
reduced time steps relative to SWE, but it has been found
that these equations can still be solved
efficiently with explicit methods, unlike SGN.
The algorithms presented in this paper and the resulting GeoClaw code can
provide useful reference solutions for the further exploration of
hyperbolic approaches.

\section*{Code and data availability}
The GeoClaw setup for the examples presented in \cref{sec:ex}
can be found in a Github repository created for this paper, 
\url{https://github.com/rjleveque/ImplicitAMR-paper},
\add{and is permanently archived on Zenodo \cite{zenodo}.}
This includes standard Clawpack-style
{\tt setrun.py} files that show the AMR resolutions, refinement regions,
and other parameters used in these examples.
\add{The software described in this paper and used for the experiments presented
here has recently been merged into GeoClaw
and is part of release v5.10.0. \cite{clawpack}}

\remove{Full reproducibility is not yet
supported since the new 1D and 2D GeoClaw code that implements
the Boussinesq equations (and implicit AMR in 2D) is still under active
development and will soon be be refactored and better documented for merging
into GeoClaw. We plan to update the paper repository in the future to make it
easier for others to reproduce the results presented here.}

\section*{Acknowledgments}
Barry Smith was instrumental in building  the new version of PETSc
(now included in v3.20)  that
allows OpenMP code in combination with MPI, and assisting us in its
use. 
We had helpful discussions with Donna Calhoun on the implementation
of the SGN equations.
St\'ephane Popinet provided information about his simulations in 
\cite{Popinet:2015} and about his Basilisk code.

This work was partially supported by the NASA Asteroid Threat Assessment
Project (ATAP) through the Planetary Defense Coordination Office and
BAERI contract AO9667.

\bibliographystyle{siamplain}
\bibliography{references}
\end{document}